\documentclass[12pt,draftcls,onecolumn]{IEEEtran}
\usepackage{cite}
\usepackage{amsmath,amssymb,amsfonts}
\usepackage{algorithmic}
\usepackage{graphicx}
\usepackage{algorithm,algorithmic}
\usepackage{mathrsfs,amssymb,amsmath,array,color}
\usepackage{bm}
\usepackage{mathbbold,bbm}
\usepackage{graphicx,epstopdf} 
\usepackage{caption}
\usepackage{graphicx,epstopdf} 
\usepackage[all]{xy}
\usepackage{chemarrow}
\makeatletter
    
    \newcommand{\Rmnum}[1]{\expandafter\@slowromancap\romannumeral #1@}

\onecolumn

\newtheorem{theorem}{Theorem}[section]
\newtheorem{definition}[theorem]{Definition} 
\newtheorem{remark}[theorem]{Remark} 
\newtheorem{lemma}[theorem]{Lemma}
\newtheorem{corollary}[theorem]{Corollary}

\newtheorem{example}{Example}

\begin{document}

\title{Stability Analysis of Biochemical Reaction Networks Linearly Conjugated to complex balanced Systems with Time Delays Added}

\author{Xiaoyu Zhang, Shibo He, Chuanhou Gao, and Denis Dochain 
\thanks{Manuscript received \today. This work was funded by the National Nature Science Foundation of China under Grant No. 62303409, 12320101001,  and 12071428, and the China Postdoctoral Science Foundation under Grant No. 2023M733115. This work is an extension of our earlier paper in the 4th IFAC Workshop on Thermodynamics Foundations of
Mathematical Systems Theory, Jul. 25-27, 2022, Canada. }
\thanks{X. Zhang and S. He are with the College of Control Science and Engineering, C. Gao is with the School of Mathematical Sciences, Zhejiang University, Hangzhou, China (e-mail:Xiaoyu\_Z@zju.edu.cn, s18he@zju.edu.cn, gaochou@zju.edu.cn (correspondence)). }

\thanks{D. Dochain is with ICTEAM, UCLouvain, B\^{a}timent Euler, avenue Georges Lema\^{i}tre 4-6, 1348 Louvain-la-Neuve, Belgium (e-mail: denis.dochain@uclouvain.be)}
}

\maketitle

\begin{abstract}
Linear conjugacy offers a new perspective to broaden the scope of stable biochemical reaction networks to the systems linearly conjugated to the well-established complex balanced mass action systems ($\ell$cCBMASs). This paper addresses the challenge posed by time delay, which can disrupt the linear conjugacy relationship and complicate stability analysis for delayed versions of $\ell$cCBMASs (D$\ell$cCBMAS).
Firstly, we develop Lyapunov functionals tailored to some D$\ell$cCBMASs by using the persisted parameter relationships under time delays. Subsequently, we redivide the phase space as several invariant sets of trajectories and further investigate the existence and uniqueness of equilibriums in each newly defined invariant set. This enables us to determine the local asymptotic stability of some D$\ell$cCBMASs within an updated framework.
Furthermore, illustrative examples are provided to demonstrate the practical implications of our approach.
\end{abstract}


\IEEEpeerreviewmaketitle
\section{Introduction}
\label{sec:introduction}
\IEEEPARstart{B}{iological} processes are inherently characterized by high non-linearity and complexity \cite{freeman2001biocomplexity}, necessitating a network framework to describe interactions among components \cite{Tang2009}. 
Biochemical Reaction Network Theory (BRNT) under this framework establishes nonlinear ordinary differential equations to describe the evolution of the biological systems and  reveals the dynamical properties from the view of network structures \cite{Feinberg2019}. Hence, BRNT falls within the domain of systems biology and serves as a fundamental pillar of synthetic biology.
Nevertheless, time delay phenomena are widespread in biological processes, functioning as critical tools for model simplification \cite{hangos2018}, regulation of biochemical processes \cite{Ri2021}, stabilization of ecosystems \cite{Li2023}, among other functions. Therefore, the exploration of delayed Biochemical Reaction Networks (BRNs) becomes essential.

Biochemical systems encompass various types of cells, biomolecules, and intricately intertwined interactions among them. Stability is a crucial property to ensure the reliability of biological systems, which is essential for the survival and proper functioning of organisms. Therefore, stability analysis is of paramount importance in BRNT \cite{Giordano2014,Fang2019, Fang2020}, as it needs to address the inherent complexity and dynamism of biological systems.
In the early stages of stability analysis, the primary focus was on establishing a connection from network topology to equilibrium distribution, followed by an exploration of equilibrium stability. Notably, there was a strong emphasis on the weakly reversible structure, a prerequisite for complex balanced mass action systems (CBMAS). Feinberg, Horn and Jackson formulated the well-known deficiency zero theorem \cite{Horn1972, Feinberg1995}, asserting that a weakly reversible deficiency zero MAS is not only complex balanced but also possesses a unique equilibrium in each positive stoichiometric compatibility class (an invariant set of the trajectory). Moreover, within the complex balanced system, each equilibrium is locally asymptotically stable, with the pseudo-Helmholtz free energy function serving as the Lyapunov function \cite{Feinberg1972}. Hangos further confirmed that complex balanced systems can maintain stability under arbitrary time delays by constructing a delayed version of the pseudo-Helmholtz free energy functional \cite{Hangos2018}. Additionally, the existence and uniqueness of the equilibrium also apply to delayed complex balanced mass action systems (DCBMASs).
However, this function or functional is limited in systems with non-weakly reversible structures. Hence, there is an urgent need for novel techniques to expand our understanding of biochemical reaction network systems' stability beyond the currently known scope.

However, analyzing a specific class of systems based on their structure is exceptionally challenging due to a lack of appropriate mathematical tools. One viable approach involves exploring the form of an unknown system by leveraging the relationship between the trajectories of a known stable system, a concept akin to system identification in contemporary biology research. System identification involves deducing interactions among biochemical molecules through experimental data collected in biochemical experiments \cite{Alonso2012, Ye2023}. Despite potentially differing structures and parameters, identified networks exhibit clear trajectories' relationships under identical initial conditions, such as dynamical equivalence and linear conjugacy. These relationships eliminate the need to explicitly consider the system's structure to infer its properties, representing a ``black-box" method for stability analysis.

For non-delayed biochemical systems, Johnston et al. \cite{Siegel2011} presented the stability of a class of biochemical reaction systems transformable into a complex balanced system ($\ell c$CBMAS) through a linear conjugacy approach where the linear transformation matrix is positively diagonal and dynamic equivalence is a special case of linear conjugacy where the corresponding matrix is a identity matrix. A natural question arises: can this linear conjugate relationship be preserved under arbitrary time delays? The answer, as detailed in our last conference paper \cite{Zhang2022}, is no. Time delays disrupt the relationship between trajectories, challenging the straightforward ``black-box" stability analysis of delayed $\ell$cCBMAS (D$\ell$cCBMAS). Further investigation is required to integrate the specific characteristics of the target system.

Despite disruptions, certain parameter relationships persist in D$\ell c$CBMAS, shedding light on their dynamic behavior. We decompose the dynamical equations, leveraging characteristics of the corresponding  complex balanced systems to establish Lyapunov functionals ensuring stability of some D$\ell c$CBMASs. However, the existence of degenerate equilibrium points raises questions about local asymptotic stability relative to stoichiometric compatibility classes (SCCs).
To address this challenge, we further re-decompose the phase space into equivalence classes, where each class represents an invariant set of system trajectories. The primary contributions of this work include:
\begin{itemize}
	\item Stability analysis of some delayed high-frequency biochemical reaction networks compared to some CBMAS to which is dynamically equivalent.
	\item Stability analysis of D$\ell c$CBMAS  whose all reactant complex are in the form of the multiple of a single species.
	\item Redivision the phase space into a set of invariant sets and derive the existence, uniqueness, and local asymptotic stability of the equilibrium relative to the newly defined invariant sets.  
\end{itemize}

Now we elucidate the overall structure of this paper. Section \ref{II} provides fundamental concepts regarding chemical reaction networks, dynamics. Section \ref{sec:3} outlines the motivation behind this paper and introduces several modules commonly used in the construction of Lyapunov functionals in the subsequent sections. Section \ref{sec:4} focuses on the discussion of Lyapunov stability for D$\ell c$CBMAS. Finally, Section \ref{sec:5} deals with a comprehensive exploration of the system's local asymptotic stability.

~\\~\\
\noindent{\textbf{Mathematical Notation:}}\\
\rule[1ex]{\columnwidth}{0.8pt}
\begin{description}
\item[\hspace{-2em}{$\mathbb{R}^n_{\geq 0},\mathbb{R}^n_{>0}:$}] $n$-dimensional non-negative real space; $n$-dimensional positive real space.
\item[\hspace{-0.6em}{$\bar{\mathscr{C}}_{+}, \mathscr{C}_{+}$}]: $\bar{\mathscr{C}}_{+}=C([-\tau,0];\mathbb{R}^{n}_{\geq 0}), \mathscr{C}_{+}=C([-\tau,0];\mathbb{R}^{n}_{>0})$ the non-negative, positive continuous function vectors defined on the interval $[-\tau,0]$, respectively.
\item[\hspace{-0.2em}{$x^{y_{\cdot i}}$}]: $x^{y_{\cdot i}}\triangleq\prod_{j=1}^{n}x_{j}^{y_{ji}}$, where $x,y_{\cdot i}\in\mathbb{R}^{n}$.
	\item[\hspace{-0.5em}{$\mathrm{Ln}(x)$}]: $\mathrm{Ln}(x)\triangleq\left(\ln{x_{1}}, \cdots, \ln{{x}_{n}} \right)^{\top}$, where $x\in\mathbb{R}^{n}_{>{0}}$.
 \item[\hspace{-0.5em}{$k^{(y,\tau)}_i$}]: If $y_{.i}=y$ and $\tau_i=\tau$, $k^{(y,\tau)}_i=k_i$. Otherwise, $k_i^{(y,\tau)}=0$.
\item[\hspace{-0.5em}{$k^{(y)}_i$}]: If $y_{.i}= y$, $k^{(y)}_i=k_i$. Otherwise, $k_i^{(y)}=0$.
\item[\hspace{-0.5em}{$n,p,r$}]: The number of species, complex, reactions of the corresponding system.
\end{description}
\rule[1ex]{\columnwidth}{0.8pt}

\section{Preliminaries}\label{II}
In this section, some basic concepts about BRN and the corresponding delayed version are reviewed \cite{Feinberg2019,Hangos2018}. 

Consider a network with $n$ species $X_1,...,X_n$, labeled by $\mathcal{S}=\{X_j, j=1,\cdots, n\}$, and $r$ reactions $R_1,...,R_r$, denoted by $\mathcal{R}=\{R_i, i=1,\cdots, r\}$, in which the $i$-th reaction follows
	$$R_i:~~~~\sum_{j=1}^n y_{ji}X_j \to \sum_{j=1}^n y'_{ji}X_j.$$ 
Here, $y_{ji}$ and $y'_{ji}$ represent the stoichiometric coefficients of species $X_j$ in the reactant and product of reaction $R_i$. The corresponding integer vector $y_{.i}\in\mathbb{Z}^n_{\geq 0}$ and $y'_{.i}\in\mathbb{Z}^n_{\geq 0}$ are referred to as the reactant and product complexes of $R_i$, respectively. All complexes form a set $\mathcal{C}=\cup_{i=1}^r\{y_{.i},y'_{.i}\}$, and all reactant complexes are denoted by $\mathcal{RC}$. The change of all species causing by the $R_i$ is depicted by its \textbf{reaction vector} $v_{.i}=y'_{.i}-y_{.i}$, all of which span a subspace of $\mathbb{R}^n$ known as the stoichiometric subspace
\begin{equation*}
	\mathscr{S}=\mathrm{span}\{y'_{\cdot i}-y_{\cdot i}:~i=1,\cdots,r\}.
\end{equation*}
We often refer to a BRN as a triple $\mathcal{N}=(\mathcal{S,C,R})$. 

Mass action kinetics is the most frequently-used law to evaluate the reaction rate, by which each reaction rate is proportional to the product of all reactant concentrations. For $R_i$, the reaction rate is $K_i(x(t))\triangleq k_ix^{y_{.i}}(t)$,
	where $x\in\mathbb{R}^n_{\geq 0}$ is the concentration vector of all species and $k_i\geq 0$ is known as the reaction rate constant.  

\begin{definition}[MAS]
An $\mathcal{N}=(\mathcal{S,C,R})$ equipped with mass action kinetics is termed a mass action system (MAS), labeled by a quadruple $\mathcal{M}=(\mathcal{S,C,R},\bm{k})$. 
\end{definition}

The dynamics of $\mathcal{M}=(\mathcal{S,C,R},\bm{k})$ thus follows $\dot{x}_{\text{nd}}(t)=f_{\text{nd}}(x)$ where 
\begin{equation}\label{eq:mas}
\begin{split}
	&f_{\text{nd}}(x)=\sum^{r}_{i=1}k_ix^{y_{\cdot i}}(y'_{\cdot i}-y_{\cdot i})\triangleq\sum_{y\in\mathcal{RC}}f^{(y)}_{\text{nd}}(x)=\sum_{y\in\mathcal{RC}}k^{(y)}_ix^{y}(y'_{\cdot i}-y),~t\geq 0.
 \end{split}
\end{equation}
Further, as the consumption of each delayed reaction occurs immediately, while the product of $R_i$ at time $t$ is actually produced at time $t-\tau_i$ ($\tau_i$ is the time delay of $R_i$), the dynamics of delayed mass action system (DeMAS, denoted by $\mathcal{M}_{\text{De}}=(\mathcal{M},\bm{\tau})$) expressed as $\dot{x}(t)=f(x)$ for any $t\geq 0$ is in the following form
\begin{equation}\label{eq:dde}
\begin{split}
&f(x)=\sum^{r}_{i=1}k_i[(x(t-\tau_i))^{y_{\cdot i}}y'_{\cdot i}-(x(t))^{y_{\cdot i}}y_{\cdot i}]\triangleq \sum_{y\in \mathcal{RC}}f^{(y)}\triangleq \sum_{y\in \mathcal{RC}}\sum_{\tau\in\bm{\tau}}f^{(y,\tau)}(x)
\end{split}
\end{equation}
where 
\begin{equation}
f^{(y,\tau)}(x)=\sum_{i=1}^r[k^{(y,\tau)}_ix^y(t-\tau)y'_{.i}-k^{(y,\tau)}_ix^y(t)y].
\end{equation}
Note that when all $\tau_i$ equal zero, \eqref{eq:dde} reduce to the non-delayed case \eqref{eq:mas}.
The phase space of $\mathcal{M}_{\text{De}}$ is given by $\bar{\mathscr{C}}_+\triangleq C([-\tau_{\text{max}}, 0]; \mathbb{R}^n_{\geq 0})$, and its interior is $\mathscr{C}_+\triangleq C([-\tau_{\text{max}}, 0]; \mathbb{R}^n_{>0})$, where $\tau_{\text{max}}=\max\{\tau_i, i=1,\cdots,r\}$.  Each trajectory of \eqref{eq:dde}, denoted as $x(t)^\psi$, always stays in some invariant set known as the stoichiometric compatibility class (SCC) once the initial condition $\psi\in \bar{\mathscr{C}}_+$ is fixed.

\begin{definition}[SCC]\label{def:SCC}
 The non-negative SCC containing  $\psi\in\bar{\mathscr{C}_{+}}$  of $\mathcal{M}_{\mathrm{De}}$ is
\begin{equation}\label{eq:Dpsi}
		\mathcal{D}_{\psi}:=\{\phi\in \bar{\mathscr{C}}_{+} \vert c_{b}(\phi)=c_{b}(\psi),~\forall ~b\in \mathscr{S}^{\bot}\},
	\end{equation} 
	where $\mathscr{S}^{\bot}$ is the orthogonal complement of $\mathscr{S}$, and the functional $c_{b}:\bar{\mathscr{C}_{+}}\rightarrow \mathbb{R}$ is 
\begin{equation}
    \begin{split}
\label{eq:f}
		c_{b}(\psi):=&b^{\top}\left[\psi(0)+\sum_{y_{\cdot i}\rightarrow y'_{\cdot i}}\left(\int_{-\tau_{i}}^{0}k_i\psi(s)^{y_{\cdot i}}ds\right)y_{\cdot i}\right]=b^{\top}g(\psi).
	\end{split}
	\end{equation}
Specially, $\mathcal{D}_\psi^+$ is a positive SCC containing $\psi$ if $\phi$ in \eqref{eq:Dpsi} is defined in $\mathscr{C}_+$.
\end{definition}

The network structure plays an important role on the dynamic behavior of a network. Thereinto, the weakly reversible structure is the most active one, and induces a class of well-known network systems, called complex balanced MASs (CBMASs).



\begin{definition}[weakly reversible]\label{def:2reversible} A BRN $\mathcal{N}$ is weakly reversible if for any reaction $y_{\cdot i}\to y'_{\cdot i}\in\mathcal{R}$, there must be a chain of reactions from $y'_{\cdot i}$ to $y_{\cdot i}$ in $\mathcal{R}$, i.e., $\exists ~y'_{\cdot i} \to y_{\cdot i_{1}}\in\mathcal{R}$, $\cdots$, $y_{\cdot i_{m}}\to y_{\cdot i}\in\mathcal{R}$ with $m\leq r$.	
\end{definition}

\begin{definition}[CB equilibrium]\label{def:equilibrium}
	An $\bar{x}\in\mathbb{R}^n_{\geq 0}$ is an equilibrium of $\mathcal{M}_{\mathrm{De}}$ if $\dot{x}\vert_{\bar{x}}=0$. It is further a CB equilibrium if
\begin{equation}\label{eq:cb}
	\sum\limits_{\{i:~y_{\cdot i}=y\}}k_i(\bar{x})^y=\sum\limits_{\{i:~y'_{\cdot i}=y\}}k_i(\bar{x})^{y_{.i}},~\mathrm{for~each~} y\in\mathcal{C}.
\end{equation}
\end{definition}

The definition emphasizes that at each CB equilibrium, the system attains a balance relative to each complex. The MAS contains a CB equilibrium is called a CBMAS, and moreover, each equilibrium in CBMAS is a CB equilibrium. While complex balanced systems indeed impose constraints on the reaction rate constants, within this category, there exists a specific subclass known as weakly reversible, zero-deficiency systems, for which the complex balancing is always true irrespective of the reaction rate constants. Here, the deficiency is defined as follows.


\begin{definition}[deficiency]
	Consider a chemical reaction network $\mathcal{N}$ with its associated graph $G=(V,E)$ where each complex is a node in $V$ and each reaction $R_i$ is an edge in $E$. Then the deficiency of $\mathcal{N}$ is defined as $p-l-s$, where $p$, $l$, and $s$ denote the number of complexes, connected components in $G$, and the dimension of the stoichiometric subspace $\mathscr{S}$, respectively.
\end{definition}
 
Horn and Jackson \cite{Horn1972} demonstrated the existence, uniqueness and further local asymptotic stability of the positive equilibrium in each SCC of any CBMAS. 
The same result has been extended to delayed CBMAS (DCBMAS) by using  $V:\mathscr{C}_{+}\rightarrow \mathbb{R}_{\geq 0}$
\begin{equation}\label{eq:Vd}
	\begin{split}
		&V(\psi)
			\triangleq \sum_{j=1}^{n}h(\psi_j(0);\bar{x}_j)+\sum_{i=1}^{r}k_i\int_{-\tau_i}^{0}h(\psi(s)^{y_{.i}};\bar{x}^{y_{.i}})ds,\\
	\end{split}
\end{equation}
where $\bar{x}$ is an equilibrium of DCBMAS and $h(\cdot;\beta):\mathbb{R}_{> 0}\to \mathbb{R}_{\geq 0}$ is defined as $h(z;\beta)=\beta-z-z\ln{\frac{\beta}{z}}$ with $z\in\mathbb{R}_{> 0}$ and the positive constant parameter $\beta$.

\section{Preparatory Work}\label{sec:3}
In this section, we elaborate on the motivation behind this work. Additionally, some Lyapunov functionals of supporting modules that play a vital role in the subsequent stability analysis are presented.
\subsection{Problem Motivation}
\label{sec:moti}
\begin{figure*}[h]
 	\centering
\includegraphics[height=6cm,width=13cm]{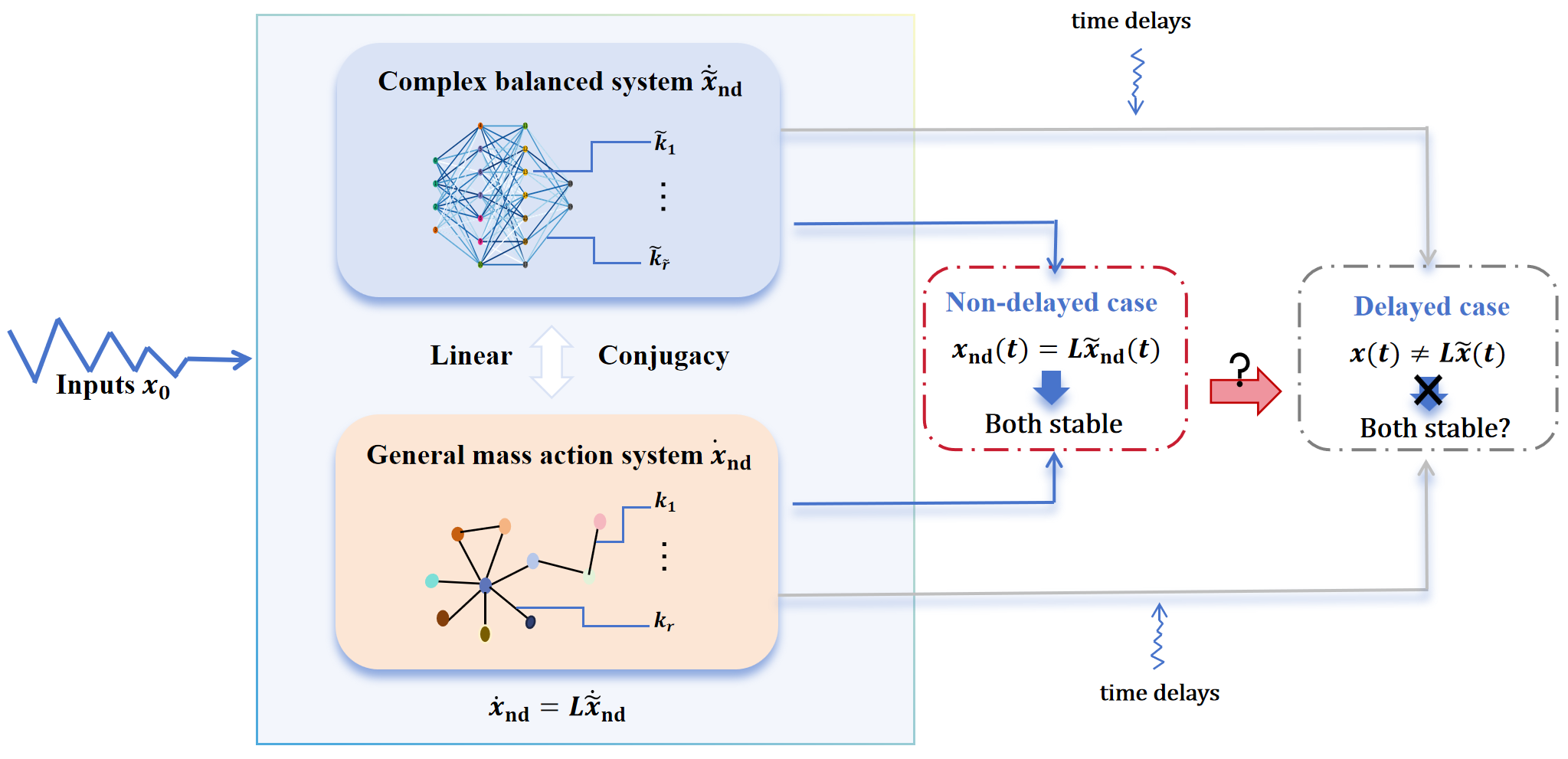}
 	\caption{The motivation graph of this article. }
 	\label{fig:mg}
 \end{figure*}

A fixed dynamics can induce various system identifications with distinct network structures and parameters. They can produce identical outputs under the same input conditions ($x_{\text{nd}}(t)=\tilde{x}_{\text{nd}}(t)$). The relationship between these systems is also referred to as dynamic equivalence. Further, dynamic equivalence can be extended to linear conjugacy, as defined below.
\begin{definition}[\cite{Siegel2011}]\label{def:31}
The two mass-action systems $\mathcal{M}$ and $\mathcal{\tilde{M}}$ are linearly conjugate if there exists a linear, bijective mapping: $\bm{\text{u}}:\mathbb{R}^{n}_{>0}\to\mathbb{R}^{n}_{>0}$ such that any trajectories $\Phi$ and $\tilde{\Phi}$ of $\mathcal{M}$,~$\tilde{\mathcal{M}}$ satisfy $\bm{\text{u}}(\Phi(x_0,t))=\tilde{\Phi}(\text{u}(x_0),t)$ for any $x_0\in\mathbb{R}^{n}_{>0}$. For non-delayed case, it also means that for each 
$y\in\mathcal{RC}\cup\mathcal{R\tilde{C}}$, there exists
\begin{equation}\label{eq:31}
\sum_{i=1}^rk_i^{(y)}(y'_{.i}-y)=\sum_{\tilde{i}=1}^{\tilde{r}}\tilde{k}_{\tilde{i}}^{(y)}\prod_{j=1}^n l_j^{-y_{j}}(\tilde{y}'_{.\tilde{i}}-y).
\end{equation}
\end{definition}
The outputs of systems under linear conjugacy satisfy $x_{\text{nd}}=L\tilde{x}_{\text{nd}}(t)$, where $L=\text{diag}(l_1,\cdots,l_n)$ is a positive, diagonal matrix. Thus, each $\mathcal{M}$ linearly conjugated to some CBMAS $\mathcal{\tilde{M}}$, known as a $\ell$cCBMAS, is stable solely based on the linear relationship between outputs, without considering the specific structure of $\mathcal{M}$. This is illustrated in the red box in Fig. \ref{fig:mg}.
Thus it is actually a ``black box" method.

 When time delays are introduced, the delayed CBMAS $\mathcal{\tilde{M}}_{\text{De}}$ remains stable according to \cite{Hangos2018}. However, the stability of the delayed $\ell$cCBMAS (D$\ell$cCBMAS), denoted as $\mathcal{M}_{\text{De}}$, cannot be directly inferred from $\tilde{\mathcal{M}}_{\text{De}}$ due to the lack of linear relationship between trajectories, namely, $\mathcal{M}_{\text{De}}$ may not be linear conjugate to some DCBMAS ($\ell$cDCBMAS). This is shown in the gray box of Fig. \ref{fig:mg} and Example \ref{ex:1}.
This indicates that solely relying on a ``black-box" method is no longer effective. Therefore, additional internal structural information within $\mathcal{M}_{\text{De}}$ needs to be incorporated to achieve stability. Thus, this paper deals with certain types of $\mathcal{M}_{\text{De}}$ by using their specific structures and the parameter relationship shown in \eqref{eq:31}.
\begin{example}\label{ex:1}
 We illustrate the motivation through a simple delayed birth-death process: 
 \begin{equation}\mathcal{M}:	2S_1\xrightarrow{k_1}S_1,~~\emptyset \xrightarrow{k_2}S_1
 \end{equation}
which is dynamically equivalent to a CBMAS 
 \begin{equation}\mathcal{\tilde{M}}: \xymatrix{2S_1\ar @{ -^{>}}^{~0.5k_1}@< 1pt> [r]& \emptyset \ar  @{ -^{>}}^{~0.5k_2}  @< 1pt> [l]}.
 \end{equation}
 The delayed version of $\mathcal{M}$ denoted as $\mathcal{M}_{\mathrm{De}}$ has the dynamics of 
\begin{equation}
\dot{x}_1=k_1x_1^2(t-\tau_1)-2k_1x_1^2(t)+k_2
\end{equation}
with $\tau_1>0$ to be the delay of the first reaction in $\mathcal{M}$, and $\mathcal{\tilde{M}}_{\mathrm{De}}$ has the dynamics of $$\dot{\tilde{x}}_1=-k_1\tilde{x}_1^2(t)+k_2$$ with any delay in reactions. These two dynamics are different for any $\tau_1\neq 0$, so the time delays disrupt the relationship of dynamical equivalence.
 	\end{example}

\subsection{Stability of some supporting modules}
Some supporting modules and their Lyapunov functionals are introduced in this subsection. 

\textbf{Supporting Module 1:}  $y\to y$ is a circular structure whose dynamics is zero for the non-delayed case. Its dynamics is entirely different in a delayed case
  \begin{equation}
  \dot{x}=kx(t-\tau)^{y}y-kx(t)^{y}y.
  \end{equation}
It is, in fact, a special and small DCBMAS itself, where the inflow and outflow of complex $y$ are always the same at each constant function $x_c\in\bar{\mathscr{C}}_+$. Thus, we can directly write the Lyapunov functional of this module using \eqref{eq:Vd}.

 %

\textbf{Supporting Module 2:} Now we define a virtual DCBMAS.
\begin{definition}[Quasi-DCBMAS]\label{def:qd} $\mathcal{\tilde{M}}_{\mathrm{q}}=(\mathcal{\tilde{M}}, \bm{\tilde{\tau}}^{\mathrm{q}})$ is called a Quasi-DCBMAS of a CBMAS $\mathcal{\tilde{M}}$ if there exists some reaction $\tilde{R}_{\tilde{i}}\in\mathcal{\tilde{R}}$ with more than one time delay. 
\end{definition}

For the delayed version of $\tilde{\mathcal{M}}$, $\mathcal{\tilde{M}}_{\text{De}}=(\mathcal{\tilde{M}},\bm{\tilde{\tau}})$, there is a maximum of $r$ distinct time delays in $\mathcal{\tilde{M}}_{\text{De}}$, precisely $\vert \bm{\tilde{\tau}}\vert \leq r$. However, the number of time delays in $\mathcal{\tilde{M}_{\text{q}}}$ denoted as $\vert \bm{\tilde{\tau}^{\text{q}}}\vert$ is not constrained by $r$. 
Further by denoting $\tilde{k}_{.\tilde{i}}^{(\tilde{\tau})}$ as the reaction rate of the $\tilde{i}$-th reaction with time delay $\tilde{\tau}$, thus 
it satisfies that $\sum_{\tilde{\tau}\in\bm{\tilde{\tau}^{\text{q}}}}\tilde{k}^{(\tilde{\tau})}_{\tilde{i}}=\tilde{k}_{\tilde{i}}$.

Then the dynamics of $\tilde{\mathcal{M}}_\mathrm{q}$ can be written as
\begin{equation}
\begin{split}
&\dot{\tilde{x}}^{\text{q}}=\sum_{\tilde{i}=1}^{\tilde{r}}\sum_{\tilde{\tau}\in \bm{\tilde{\tau}}^{\mathrm{q}}}\tilde{k}^{(\tilde{\tau})}_{\tilde{i}}
\left[(x(t-\tilde{\tau}))^{\tilde{y}_{\cdot \tilde{i}}}\tilde{y}'_{\cdot \tilde{i}}-(x(t))^{\tilde{y}_{\cdot \tilde{i}}}\tilde{y}_{\cdot \tilde{i}}\right]\\
&=\sum_{\tilde{y}\in \mathcal{R\tilde{{C}}}}\sum_{\tilde{\tau}\in \bm{\tilde{\tau}}^{\mathrm{q}}}\sum_{\tilde{i}=1}^{\tilde{r}}\tilde{k}^{(\tilde{y},\tilde{\tau})}_{\tilde{i}}
\left[(x(t-\tilde{\tau}))^{\tilde{y}}\tilde{y}'_{\cdot \tilde{i}}-(x(t))^{\tilde{y}}\tilde{y}\right]
\end{split}
\end{equation}
It is easy to verify that the following $V_{\text{q}}$ is the Lyapunov functional of $\mathcal{\tilde{M}}_{\text{q}}$ 
\begin{equation}\label{eq:LCBVd}
		\begin{split}
			V_{\text{q}}(\psi)=\sum_{j=1}^{\tilde{n}}h(\psi_j(0); \bar{x}_j)+\sum_{\tilde{y}\in \mathcal{R\tilde{C}}}\sum_{\tilde{\tau}\in\bm{\tilde{\tau}^{\text{q}}}}\sum_{\tilde{i}=1}^{\tilde{r}}\tilde{k}_{\tilde{i}}^{(\tilde{y},\tilde{\tau})}\int^{0}_{-\tilde{\tau}}h(\psi(s)^{\tilde{y}};\bar{x}^{\tilde{y}}) ds.
		\end{split}
	\end{equation}
Here,  $\tilde{k}_{\tilde{i}}^{(\tilde{y},\tilde{\tau})}$ is equal  to $\tilde{k}^{(\tilde{\tau}) }_{\tilde{i}}$ if $\tilde{y}=\tilde{y}_{\cdot\tilde{i}}$, otherwise, it's zero. The equality of $\dot{V}_{\text{q}}(t)\leq 0$ is achieved iff $x$ is an equilibrium of $\tilde{\mathcal{M}}_{\text{q}}$, which is also an equilibrium of $\mathcal{\tilde{M}}$.

\textbf{Supporting module 3:} For the above mentioned $\mathcal{M}_{\text{De}}$ and $\mathcal{\tilde{M}}_{\text{De}}$, if $\mathcal{M}_{\text{De}}$ is linearly conjugate to $\mathcal{\tilde{M}}_{\text{De}}$, called $\ell c$DCBMAS, the Lyapunov-stability of $\mathcal{M}_{\text{De}}$ was reached \cite{Zhang2022}. Here, we generalize this to $\ell$cQuasi-DCBMAS, denoted as $\mathcal{M}_{\text{ql}}$, which is linearly conjugate to a Quasi-DCBMAS $\mathcal{\tilde{M}}_{\text{q}}$. Namely, the dynamics of $\mathcal{M}_{\text{ql}}$ and $\mathcal{\tilde{M}}_{\text{q}}$, denoted as $\dot{x}^{\text{ql}}$ and $\dot{\tilde{x}}^{\text{q}}$, respectively, satisfy the relation $\dot{x}^{\text{ql}}=L\dot{\tilde{x}}^{\text{q}}.$ 
\begin{lemma}[stability of $\ell$cQuasi-DCBMAS]
If $\mathcal{M}_{\mathrm{ql}}$ is an $\ell$cQuasi-DCBMAS, i.e., being linearly conjugate to a Quasi-DCBMAS $\mathcal{\tilde{M}}_{\mathrm{q}}$, then $\mathcal{M}_{\mathrm{ql}}$ is Lyapunov stable.
\end{lemma}  
\textit{\textbf{Proof:}}
 It is easy to verify the Lyapunov-stability of $\mathcal{M}_{\text{ql}}$ by using the same process of the stability analysis of $\ell$cDCBMAS shown in Theorem 5 in \cite{Zhang2022} with the following $V_{\text{ql}}:\mathscr{C}_+\to \mathbb{R}_{\geq 0}$ serving as its Lyapunov functional.
\begin{equation*}
\begin{split}
		V_{\text{ql}}(\psi)&=\sum_{j=1}^{n}l_j^{-1}h(\psi_j(0); \bar{x}_j)
+\sum_{\tilde{y}\in \mathcal{R\tilde{C}}}\sum_{\tilde{i}=1}^{\tilde{r}}\sum_{\tilde{\tau}\in\bm{\tilde{\tau}^{\text{q}}}}\tilde{k}_{\tilde{i}}^{(\tilde{y},\tilde{\tau})}\prod_{j=1}^{n}l_j^{-\tilde{y}_{j}}\int^{0}_{-\tilde{\tau}}h(\psi(s)^{\tilde{y}};\bar{x}^{\tilde{y}}) ds.\end{split}
\end{equation*}
	$\hfill\blacksquare$
 
\section{Stability analysis of some D$\ell$cCBMASs}\label{sec:4}
In this section, we utilize those three supporting modules to study the stability of some D$\ell$cCBMASs, i.e., $\mathcal{M}_{\text{De}}=(\mathcal{M},\bm{\tau})$ with $\mathcal{M}$ to be linearly conjugate/dynamically equivalent (special linear conjugacy) to a CBMAS $\mathcal{\tilde{M}}$. Note that the $\ell c$DCBMAS class (a subset of D$\ell$cCBMASs) is naturally excluded, whose stability has been addressed in our earlier work \cite{Zhang2022}. We define three classes of D$\ell$cCBMASs to carry out analysis. 

\subsection{Stability of D$\ell$cCBMASs where the same reactant complexes have the same time delay}\label{subsec:4.1}
This part considers a class of D$\ell$cCBMASs with every $\mathcal{M}_{\text{De}}=(\mathcal{M},\bm{\tau})$ to satisfy that $\mathcal{M}$ is dynamically equivalent to a CBMAS $\tilde{\mathcal{M}}$ and the same reactant complex $y$ shares the same time delay $\tau^{(y)}$, labeled by $\mathcal{M}_{\text{De1}}$ for identification. 

\begin{theorem}\label{thm:st}
Given an $\mathcal{M}_{\mathrm{De1}}$ with an equilibrium $\bar{x}\in \mathbb{R}_{>0}^n$, if for each $y\in \mathcal{RC}$ we have $K^{(y)}\triangleq \sum_{i=1}^rk^{(y)}_i-\sum_{\tilde{i}=1}^{\tilde{r}}\tilde{k}^{(y)}_{\tilde{i}}\geq 0$, then $\bar{x}$ is stable.	
\end{theorem}
\textit{\textbf{Proof:}}
The dynamics of $\mathcal{M}_{\mathrm{De1}}$ can be written as 
\begin{equation*}
	\begin{split}
\dot{x}&=f(x)=\sum_{y\in\mathcal{RC}}\sum_{i=1}^{r}\left[k_i^{(y)}x^y(t-\tau^{(y)})y'_{.i}-k^{(y)}_ix^y(t)y\right]\\
	&=\sum_{y\in\mathcal{RC}}x^y(t-\tau^{(y)})\sum_{\tilde{i}=1}^{\tilde{r}}\tilde{k}^{(y)}_{\tilde{i}}(\tilde{y}'_{.\tilde{i}}-y)+\sum_{y\in\mathcal{RC}}\sum_{i=1}^{r}k^{(y)}_i(x^y(t-\tau^{(y)})-x^y(t))y\\
&=\underbrace{\sum_{y\in\mathcal{RC}}\sum_{\tilde{i}=1}^{\tilde{r}}\left[\tilde{k}_{\tilde{i}}^{(y)}x^y(t-\tau^{(y)})\tilde{y}'_{.\tilde{i}}-\tilde{k}^{(y)}_{\tilde{i}}x^y(t)y\right]}_{f_1}+\underbrace{\sum_{y\in\mathcal{RC}}K^{(y)}(x^y(t-\tau^{(y)})-x^y(t))y}_{f_2}.
\end{split}
\end{equation*}
Clearly, $f_1$ is exactly the dynamics of the delayed version of $\tilde{\mathcal{M}}$, $\mathcal{\tilde{M}}_{\text{De}}=\{\mathcal{\tilde{M}},\bm{\tilde{\tau}}\}$, and
$f_2$ is the sum of the dynamics of several \textbf{Supporting modules 1}. Then we construct the following functional $V(\psi) = V_1(\psi) + V_2(\psi) - \sum_{j=1}^nh(\psi_j(0))$ where $V_1(\psi)$ and $V_2(\psi)$ are the Lyapunov functionals of $f_1$ and $f_2$. And further consider its Lie derivative along each trajectory $x(t)$ of $\mathcal{M}_{\text{De}}$. It is easy to verify that $L_{f}V = L_{f_1}V_1 + L_{f_2}V_2\leq 0$ from $L_{f_i}V_i\leq 0$ for $i=1, 2$. The equality $L_{f}V(\psi)=0$ holds if and only if each $L_{f_i}V_i(\psi)$ equals zero. For $i=1$, it implies $\psi$ is a constant equilibrium; $L_{f_2}V_2(\psi)=0$ if and only if $\psi$ is a constant function on $\mathcal{\tilde{D}}_{\psi}$. Therefore, $L_{f}V$ equals zero if and only if $\psi$ is an equilibrium of $\mathcal{M}_{\text{De}}$. Consequently, each equilibrium is Lyapunov stable.
$\hfill\blacksquare$
    \begin{remark}
Theorem \ref{thm:st} also ensures the holding of the case that the number of reactant complex in $\mathcal{M}_{\text{De}}$ larger than that of $\mathcal{\tilde{M}}$. For example, if some $y_1$ is in $\mathcal{RC}$ but not in $\mathcal{R\tilde{C}}$, there exists $\sum_{\tilde{i}=1}^{\tilde{r}}\tilde{k}_{\tilde{i}}^{(y_1)}(\tilde{y}'_{\cdot \tilde{i}}-y_1)=0$. In this case, $f^{(y_1)}$ can be included into $f_{2}$. 
\end{remark}
    
    \begin{remark}\label{rem:42}
The condition on $K^{(y)}$ reflects the structural characteristics of $\mathcal{M}$. It implies that while two systems may share identical inputs (initial points) and outputs (trajectories), the actual pathways going from inputs to outputs may vary significantly. 
Note that the concentration change induced by each reaction per unit time is 
$k_ix^y(t)v^{(y)}_{.i}$. Thus when  
	$v^{(y)}_{.i}$ in $\mathcal{M}$ are significantly smaller than $\tilde{v}^{(y)}_{.\tilde{i}}$
in $\tilde{\mathcal{M}}$, attaining the same trajectory requires a higher rate of reactions in $\mathcal{M}$, namely, $K^{(y)}\geq 0$.
\end{remark}

 It is straightforward to get the Lyapunov stability of $\mathcal{M}_{\text{De}}$ in Example \ref{ex:1} by Theorem \ref{thm:st} as its $K^{(y)}$ satisfy $K^{(2S_1)}=k_1/2>0$ and $K^{(\emptyset)}=k_2/2>0$.
\subsection{Stability of some high-frequency D$\ell$cCBMASs}\label{subsec:4.2}
Inspired by Remark \ref{rem:42}, this subsection rigorously defines a class of D$\ell$cCBMASs called high-frequency DeMASs denoted as $\mathcal{M}_{\text{De2}}=(\mathcal{M},\bm{\tau})$ and analyze their stability. 
\begin{definition}\label{def:hf}
$\mathcal{M}_{\mathrm{De}}=(\mathcal{M}, \bm{\tau})$ is a HF-DeMAS denoted as $\mathcal{M}_{\mathrm{De2}}$ with respect to $\mathcal{\tilde{M}}$ if $\mathcal{M}$ is dynamically equivalent to $\mathcal{\tilde{M}}$, and the 1-norm of any weighted average vector of non-zero $v_{i}^{(y,\tau)}$ of $\mathcal{M}_{\mathrm{De}}$ with respect to reaction rate constants, denoted as $$\bar{v}^{(y,\tau)}=\sum_{i=1}^rk_i^{(y,\tau)}v_{.i}^{(y,\tau)} / \sum_{i=1}^rk_i^{(y,\tau)}$$ is no greater than the minimum value among the 1-norms of all nonzero $\tilde{v}_{\tilde{i}}^{(y)}$ of $\tilde{\mathcal{M}}$.
\end{definition}

\begin{example}[HF-DeMAS]\label{ex:hf}
Consider $\mathcal{M}_{\mathrm{De}}=(\mathcal{M}, \bm{\tau})$ and $\mathcal{\tilde{M}}$ shown in Figure \ref{fig:HF}
\begin{figure}[H]
\centering
\includegraphics[height=5cm,width=13cm]{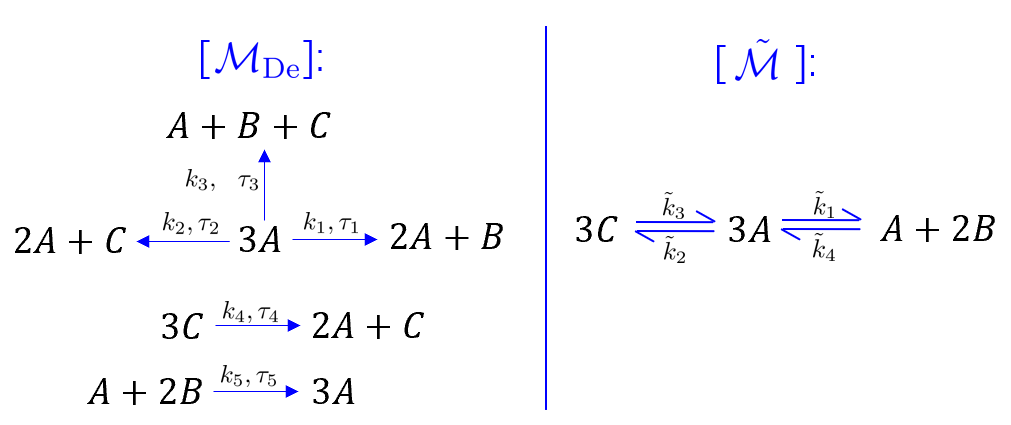}
\caption{An example of HF-DeMAS.}
\label{fig:HF}
	\end{figure}
\noindent Note that $\mathcal{\tilde{M}}$ is a weakly reversible system with zero deficiency, making it a CBMAS independent of coefficients. When $ k_i=\tilde{k}_1=\tilde{k}_4=k$, $\tilde{k}_2=\tilde{k}_3=2k/3$, the dynamics of $\tilde{\mathcal{M}}$ and the non-delayed version of $\mathcal{M}_{\mathrm{De}}$ denoted as $\mathcal{M}$ are dynamically equivalent.
	
	The reactant complex set is $\mathcal{RC}=\{3A, 3C, A+2B\}$. 
		The set of all reaction vectors with reactant complex $3A$ in $\mathcal{M}_{\mathrm{De}}$ and $\mathcal{\tilde{M}}$ are $\{(-1,0,1),(-1,1,0),(-2,1,1)\}$ and $\{(-3,0,3),(-2,2,0)\}$ respectively. If  $\tau_1=\tau_3=\tau$, the corresponding weighted average is $\bar{v}^{(y,\tau)}=(-3/2, 1, 1/2)$. The 1-norm of $\bar{v}^{(y,\tau)}$ is smaller than that of any other reaction vector $\tilde{v}^{(y)}$. If $\tau_1\neq \tau_2\neq \tau_3$, it is easy to verify that the 1-norm of each $v^{(y,\tau_i)}_i$ is no larger than the 1-norm of $\tilde{v}^{(y)}_{\tilde{i}}$. The same for other reactant complex. Thus, we can conclude that $\mathcal{M}_{\mathrm{De}}$ is a HF-DeMAS $\mathcal{M}_{\mathrm{De2}}$ compared to $\mathcal{\tilde{M}}$.
\end{example}
To investigate the stability of $\mathcal{M}_{\text{De2}}$, 
 we need to decompose its dynamics into the sum of dynamics of several supporting modules.
First, we consider the decomposition method  for $\mathcal{M}_{\mathrm{De2}}$ whose dimension of the subspace spanned by all reaction vectors $v_{i}^{(y)}$ with reactant complex $y$, denoted as $\mathscr{S}^{(y)}=\text{span}\{v^{(y)}_{.i}, i=1,\cdots,r\}$, is 1-dimensional with $w^{(y)}$ as its basis. Thus each product complex $y'_{.i}$ in $\mathcal{M}_{\text{De2}}$ and $\tilde{y}'_{.i}$ in $\mathcal{\tilde{M}}$ can be written as $y'_{.i}=y+a_i^{(y)}w^{(y)}$ and $y'_{.\tilde{i}}=y+\tilde{a}_{\tilde{i}}^{(y)}w^{(y)}$, respectively.
 
\begin{theorem}\label{thm:bz}
Consider $\mathcal{M}_{\mathrm{De2}}$ with 1d $\mathscr{S}^{(y)}$ for each $y\in \mathcal{RC}$ relative to some CBMAS $\mathcal{\tilde{M}}$. Then $\mathcal{M}_{\mathrm{De2}}$ is stable if one of the following conditions is satisfied
	\begin{itemize}
		\item each $Z^{(y,\tau)}\triangleq\sum_{i=1}^rk_i^{(y,\tau)}a_i^{(y,\tau)}$ has the same sign as $Z^{(y)}\triangleq\sum_{i=1}^rk_i^{(y)}a_i^{(y)}$.
		\item Otherwise, $\mathcal{\tilde{M}}$ is a weakly reversible, zero deficiency system where $\tilde{a}^{(y)}_{\tilde{i}_1}>0$ and $\tilde{a}^{(y)}_{\tilde{i}_2}<0$ both exist for some reactions $\tilde{R}_{\tilde{i}_1}$ and $\tilde{R}_{\tilde{i}_2}$.
	\end{itemize}
\end{theorem}
\textit{\textbf{Proof:}}
The proof of this is shown in Appendix B and the lemma severing as its technical tool are presented and proved in Appendix A.


Now we can generalize the above theorem.
\begin{corollary}\label{cor:hfc}
Consider a $\mathcal{M}_{\mathrm{De2}}$ concerning a CBMAS $\mathcal{\tilde{M}}$. And each $Z^{(y,\tau)}\triangleq\sum_{i=1}^rk_i^{(y)}v^{(y,\tau)}_{i}$ resides within the cone formed by the reaction vectors $\tilde{v}^{(y)}_{\tilde{i}}, \tilde{i}=1,\cdots, \tilde{r}$ in $\mathcal{\tilde{M}}$. Then $\mathcal{M}_{\mathrm{De2}}$ is stable if one of the following conditions is met:
	\begin{itemize}
		\item All reaction vectors $\tilde{v}^{(y)}_{\tilde{i}}$ are mutually independent.
		\item $\mathcal{\tilde{M}}$ is a weakly reversible system with zero deficiency.
	\end{itemize}
\end{corollary}
\textit{\textbf{Proof:}}
\textbf{Case I:} Each $Z^{(y,\tau)}$ can be expressed as $Z^{(y,\tau)}=\sum_{\tilde{i}=1}^{\tilde{r}}e_{\tilde{i}}^{(y,\tau)}\tilde{v}^{(y)}_{.\tilde{i}}$ where each $e_{\tilde{i}}^{(y,\tau)}\geq 0$. The dynamic equivalence between $\mathcal{M}$ and $\mathcal{\tilde{M}}$ ensures that
\begin{equation}
	\begin{split}
&Z^{(y)}=\sum_{\tau\in\bm{\tau}}Z^{(y,\tau)}=\sum_{\tau\in\bm{\tau}}\sum_{\tilde{i}=1}^{\tilde{r}}e^{(y,\tau)}_{\tilde{i}}\tilde{v}^{(y)}_{.\tilde{i}}=\sum_{\tilde{i}=1}^{\tilde{r}}e^{(y)}_{\tilde{i}}\tilde{v}^{(y)}_{.\tilde{i}}=\tilde{Z}^{(y)}=\sum_{\tilde{i}=1}^{\tilde{r}}\tilde{k}^{(y)}_{\tilde{i}}\tilde{v}^{(y)}_{.\tilde{i}}.
 \end{split}
\end{equation}
Since all $\tilde{v}^{(y)}_{.\tilde{i}}$ are independent, we can deduce that $e_{\tilde{i}}^{(y)}=\tilde{k}_{\tilde{i}}^{(y)}$ for each $\tilde{i}$. Thus,  similar to the Case I in Theorem \ref{thm:bz}, the Lyapunov functional of $\mathcal{M}_{\text{De2}}$ can be constructed. 

\textbf{Case II:}
$\mathcal{\tilde{M}}$ is a special case of the case I. Then we can derive the result directly from Case I.
$\hfill\blacksquare$
 \begin{example}
 The dynamics of $\mathcal{M}_{\mathrm{De}}$ defined in Example \ref{ex:hf} is in the following form:
 \begin{equation*}
 	\begin{split}
		\dot{x}&=k_1x_A^3(t-\tau_1)\left(\begin{matrix}
 			2\\1\\0
 		\end{matrix}\right)+k_2x_A^3(t-\tau_2)\left(\begin{matrix}
 			2\\0\\1
 		\end{matrix}\right)+k_3x_A^3(t-\tau_3)\left(\begin{matrix}
 			1\\1\\1
 		\end{matrix}\right)-3kx_A^3(t)\left(\begin{matrix}
 			3\\0\\0
 		\end{matrix}\right)\\
 		&+k_4x_C^3(t-\tau_4)\left(\begin{matrix}
 			2\\0\\1
 		\end{matrix}\right)-k_4x_C^3(t)\left(\begin{matrix}
 		0\\0\\3
 		\end{matrix}\right)+k_5x_Ax_B^2(t-\tau_5)\left(\begin{matrix}
 		3\\0\\0
 		\end{matrix}\right)-k_5x_Ax_B^2(t)\left(\begin{matrix}
 		1\\2\\0
 		\end{matrix}\right)
   \end{split}
   \end{equation*}
   \begin{equation*}
       \begin{split}
 		&=\frac{k}{2}x_A^3(t-\tau_1)\left(\begin{matrix}
 			1\\2\\0
 		\end{matrix}\right)-\frac{k}{2}x_A^3(t)\left(\begin{matrix}
 		3\\0\\0
 		\end{matrix}\right)+\frac{k}{3}x^3_A(t-\tau_2)\left(\begin{matrix}
 		0\\0\\3
 		\end{matrix}\right)-\frac{k}{3}x^3_A(t)\left(\begin{matrix}
 		3\\0\\0
 		\end{matrix}\right)\\
 		&+\frac{k}{3}x^3_A(t-\tau_3)\left(\begin{matrix}
 			0\\0\\3
 		\end{matrix}\right)-\frac{k}{3}x_A^3(t)\left(\begin{matrix}
 		3\\0\\0
 		\end{matrix}\right)+\frac{k}{2}x^3_A(t-\tau_3)\left(\begin{matrix}
 		1\\2\\0
 		\end{matrix}\right)-\frac{k}{2}x_A^3(t)\left(\begin{matrix}
 		3\\0\\0
 		\end{matrix}\right)\\
 	&+\frac{2k}{3}x^3_C(t-\tau_4)\left(\begin{matrix}
 			3\\0\\0
 		\end{matrix}\right)-\frac{2k}{3}x^3_C(t)\left(\begin{matrix}
 		0\\0\\3
 		\end{matrix}\right)+kx_Ax_B^2(t-\tau_4)\left(\begin{matrix}
 		3\\0\\0
 		\end{matrix}\right)-kx_Ax_B^2(t)\left(\begin{matrix}
 		1\\2\\0
 		\end{matrix}\right)\\
 		&+\frac{k}{2}x_A^3(t-\tau_1)\left(\begin{matrix}
 			3\\0\\0
 		\end{matrix}\right)-\frac{k}{2}x^3_A(t)\left(\begin{matrix}
 		3\\0\\0
 		\end{matrix}\right)+\frac{2k}{3}x^3_A(t-\tau_2)\left(\begin{matrix}
 		3\\0\\0
 		\end{matrix}\right)-\frac{2k}{3}x^3_A(t)\left(\begin{matrix}
 		3\\0\\0
 		\end{matrix}\right)\\
 		&+\frac{k}{6}x_A^3(t-\tau_3)\left(\begin{matrix}
 			3\\0\\0
 		\end{matrix}\right)-\frac{k}{6}x_A^3(t)\left(\begin{matrix}
 		3\\0\\0
 		\end{matrix}\right)+\frac{k}{3}x_C^3(t-\tau_4)\left(\begin{matrix}
 		0\\0\\3
 		\end{matrix}\right)-\frac{k}{3}x_C^3(t)\left(\begin{matrix}
 		0\\0\\3
 		\end{matrix}\right)
 		 	\end{split}
 \end{equation*}
 
 The first three lines of the last part of the above equation represent a Quasi-DCBMAS denoted as $\mathcal{\tilde{M}}_{\mathrm{q}}$ defined in Section \ref{sec:3} relative to $\mathcal{\tilde{M}}$ shown in Example \ref{ex:hf}. This is evident from the fact that $\sum_{\tau\in\bm{\tau}}\tilde{k}_{\tilde{i}}^{(y,\tau)}=\tilde{k}_{\tilde{i}}$, where specific values are assigned such as $\tilde{k}_1^{(3A,\tau_1)}=\frac{k}{2}$, $\tilde{k}_1^{(3A,\tau_3)}=\frac{k}{2}$, $\tilde{k}_2^{(3A, \tau_2)}=\frac{k}{3}$, $\tilde{k}_2^{(3A, \tau_3)}=\frac{k}{3}$, $\tilde{k}_3^{(A+2B, \tau_4)}=\frac{2k}{3}$, $\tilde{k}_4^{(3C,\tau_5)}=k$, and other $\tilde{k}_{\tilde{i}}^{(y,\tau_i)}=0$. Then it is Lyapunov stable from Corollary \ref{cor:hfc} and its Lyapunov functional can be derived from the linear combination of those of Supporting modules 1 and 2.
 

  \end{example}
      
\subsection{The stability of D$\ell$cCBMASs with each $y\in\mathcal{RC}$ a multiple of some single species}\label{subsec:4.3}
This part focuses on another kind of D$\ell$cCBMAS with $L=\text{diag}(l_1,\cdots,l_n)$ denoted as $\mathcal{M}_{\text{De3}}$ where each reactant complex $y\in \mathcal{RC}$ is a multiple of single species, namely, $y=y_{j_1}e_{j_1}$ for some species $X_{j_1}$ where $e_{j_1}\in \mathbb{R}^{n}_{\geq 0}$ is the indicator vector.
\begin{theorem}\label{thm:ss}
	Each $\mathcal{M}_{\mathrm{De3}}=(\mathcal{M},\tau)$ relative to some CBMAS $\mathcal{\tilde{M}}$ is Lyapunov stable.
\end{theorem}

\textit{\textbf{Proof:}}
$\mathcal{M}$ is linearly conjugate to some CBMAS $\mathcal{\tilde{M}}$ with $L=\text{diag}(l_1,\cdots,l_n)$.
As $\mathcal{\tilde{C}}$ is equal to $\mathcal{R}\mathcal{\tilde{C}}$ for $\mathcal{\tilde{M}}$, each $\tilde{y}_{.\tilde{i}}$ and $\tilde{y}'_{.\tilde{i}}$ in it are in the form of $\tilde{y}_{.\tilde{i}}=\tilde{y}_{j_1\tilde{i}}e_{j_1}$ and $\tilde{y}'_{.\tilde{i}}=\tilde{y}'_{j_2\tilde{i}}e_{j_2}$ for some species $X_{j_1}$ and $X_{j_2}$.
By letting  $\bar{k}^{(y)}_{\tilde{i}}\triangleq \tilde{k}_{\tilde{i}}^{(y)}\prod_{j=1}^n l_j^{-y_{j}}$, linear conjugacy implies that for each reactant complex $y$ 
\begin{equation*}\label{eq:42}
		\begin{split}
			&Y^{(y)}\triangleq\sum_{j=1}^nY^{(y)}_je_j=\sum_{\tilde{i}=1}^{\tilde{r}}\bar{k}_{\tilde{i}}^{(y)}L(\tilde{y}'_{.\tilde{i}}-y)+\sum_{i=1}^rk_i^{(y)}y\\
   &=\sum_{j=1}^n\sum_{\tilde{y}'_{j\tilde{i}}>0}\bar{k}^{(y)}_{\tilde{i}}l_{j}\tilde{y}'_{j\tilde{i}}e_j+\underbrace{(\sum_{i=1}^rk_i^{(y)}-\sum_{\tilde{i}=1}^{\tilde{r}}\bar{k}_{\tilde{i}}^{(y)}l_{j_1})}_{K^{(y)}}y_{j_1}e_{j_1}
		\end{split}
	\end{equation*}
Further define $Y^{(y,\tau)}\triangleq \sum_{i=1}^rk_i^{(y,\tau)}y'_{.i}=\sum_{j=1}^nY^{(y,\tau)}_je_j$. Note that there is no case where $Y^{(y,\tau)}_j\neq 0$ but $Y^{(y)}_j=0$, as their definitions reveal that they are all non-negative. Furthermore, if $Y^{(y)}_j\neq 0$, there must exist some $\tau$ such that $Y^{(y,\tau)}_j\neq 0$. 
Combining that $\sum_{\tau\in \bm{\tau}}Y^{(y,\tau)}=Y^{(y)}$, we denote $\delta_{j}^{(y,\tau)}=\frac{Y_j^{(y,\tau)}}{Y_j^{(y)}}$ for each $Y_j^{(y)}\neq 0$ which means that $\sum_{\tau\in\bm{\tau}}\delta_j^{(y,\tau)}=1$ for each $j$. 
Then $f^{(y,\tau)}$ can be written as 
	\begin{equation*}
		\begin{split}
&f^{(y,\tau)}=x^y(t-\tau)Y^{(y,\tau)}-\sum_{i=1}^rk_{i}^{(y,\tau)}x^y(t)y\\
&=\left[\sum_{j=1}^n\sum_{\tilde{y}'_{j\tilde{i}>0}}\delta_j^{(y,\tau)}\bar{k}^{(y)}_{\tilde{i}}l_j\tilde{y}'_{j\tilde{i}}e_j+\delta_{j_1}^{(y,\tau)}K^{(y)}y\right]x^y(t-\tau)-\sum_{i=1}^rk_{i}^{(y,\tau)}x^y(t)y\\&=\sum_{\tilde{i}=1}^{\tilde{r}}\bar{k}_{\tilde{i}}^{(y,\tau)}x^y(t-\tau)L\tilde{y}'_{.i}+K^{(y,\tau)}x^y(t-\tau)y-\sum_{i=1}^rk_{i}^{(y,\tau)}x^y(t)y\\
		\end{split}
	\end{equation*}
	where $\bar{k}^{(y,\tau)}_{\tilde{i}}=\sum_{j=1}^n\mathbbold{1}_{\tilde{y}'_{j\tilde{i}}\neq 0}\delta_j^{(y,\tau)}\bar{k}_{\tilde{i}}^{(y)}$ and $K^{(y,\tau)}=\delta_{j_1}^{(y,\tau)}(\sum_{i=1}^rk_i^{(y)}-\sum_{\tilde{i}=1}^{\tilde{r}}\bar{k}_{\tilde{i}}^{(y)}l_{j_1})$.
	
Then each $f^{(y)}$ is one of the four cases:

		\textbf{Case I:} If there does not exist the reaction $y_{j_1}X_{j_{1}}\to \tilde{y}'_{j_1\tilde{i}}X_{j_{1}}$ in $\mathcal{\tilde{M}}$, namely, $\sum_{\tilde{y}'_{j_1\tilde{i}}>0}\bar{k}^{(y)}_{\tilde{i}}l_{{j}_1}\tilde{y}'_{j_1\tilde{i}}e_{j_1}=0$. Further the non-negative property of each $Y_j^{(y)}$ ensures that $K^{(y)}=\sum_{i=1}^rk_i^{(y)}-\sum_{\tilde{i}=1}^{\tilde{r}}\bar{k}_{\tilde{i}}^{(y)}l_{j_1}\geq 0$. $f^{(y)}$ is exactly
	\begin{equation*}
	\begin{split}\small
&f^{(y)}=\underbrace{\sum_{\tau\in\bm{\tau}}\left[\sum_{\tilde{i}=1}^{\tilde{r}}(\bar{k}_{\tilde{i}}^{(y,\tau)}x^y(t-\tau)L\tilde{y}'_{.i}-\bar{k}_{\tilde{i}}^{(y,\tau)}x^y(t)Ly)\right]}_{f_4}+\underbrace{\sum_{\tau\in\bm{\tau}}K^{(y,\tau)}(x^y(t-\tau)-x^y(t))y}_{f_5}\\
\end{split}
\end{equation*}
\begin{equation*}
    \begin{split}
+&\underbrace{\sum_{\tau\in\bm{\tau}}\left[\sum_{\tilde{i}=1}^{\tilde{r}}\bar{k}_{\tilde{i}}^{(y,\tau)}x^y(t)Ly+K^{(y,\tau)}x^y(t)y-\sum_{i=1}^rk^{(y,\tau)}x^y(t)y\right]}_{f_6}\\
		\end{split}
	\end{equation*}
	The definition of $\bar{k}^{(y,\tau)}_{\tilde{i}}$ ensures that $f_4$ is the dynamics of reactions with reactant complex $y$ of some $\ell$cQuasi-DCBMAS denoted as $\mathcal{M}_{ql}$ of $\mathcal{\tilde{M}}$. The Lyapunov functional of $f_5$ can be derived by the supporting module 1 as each $K^{y,\tau}>0$.
	$f_6$ is zero from the form of $K^{(y,\tau)}$.
	
	\textbf{Case II:} There exist reactions in the form of $y_{j_1}X_{j_1} \to \tilde{y}'_{j_1\tilde{i}}X_{j_1}$ in $\tilde{\mathcal{M}}$. If $K^{(y)}>0$, the same proof in Case I can be used to address this situation. Now we further consider that $K^{(y)}<0$ and $Y^{(y)}_{j_1}=0$, namely, $
\sum_{\tilde{y}'_{j_1 \tilde{i}}>0}\bar{k}_{\tilde{i}}^{(y)}l_{j_1}\tilde{y}'_{j_1\tilde{i}}e_{j_1}+K^{(y)}y=0.$
 In this case, $$\sum_{\tau\in\bm{\tau}}Y^{(y,\tau)}x^y(t-\tau)=\sum_{\tau\in\bm{\tau}}\sum_{\tilde{y}'_{j_1\tilde{i}}=0}\bar{k}_{\tilde{i}}^{(y,\tau)}x^y(t-\tau)L\tilde{y}'_{.\tilde{i}}.$$
Then $f^{(y)}$ is equal to $f^{(y)}=f_7+f_8$ where
\begin{equation*}\small
f_7=\sum_{\tau\in\bm{\tau}}Y^{(y,\tau)}x^y(t-\tau)+\sum_{\tilde{y}'_{j_1{\tilde{i}}>0}}\bar{k}_{\tilde{i}}^{(y)}x^y(t)L\tilde{y}'_{.i}-\sum_{\tilde{i}=1}^{\tilde{r}}\bar{k}_{\tilde{i}}^{(y)}x^y(t)Ly
\end{equation*}
and 
$$f_8=-\sum_{\tilde{y}'_{j_1{\tilde{i}}>0}}\bar{k}_{\tilde{i}}^{(y)}x^y(t)L\tilde{y}'_{.\tilde{i}}+(\sum_{\tilde{i}=1}^{\tilde{r}}\bar{k}_{\tilde{i}}^{(y)}l_{j_1}-\sum_{i=1}^rk_i^{(y)})x^y(t)y$$
$f_7$ is the dynamics of reactions with reactant complex $y$ of some $\ell$cQuasi-DCBMAS denoted as $\mathcal{M}_{\mathrm{ql}}$ defined in Section \ref{sec:3} of $\mathcal{\tilde{M}}$. The last line of this equation denoted as $f_8$ is exactly zero from the definition $K^{(y)}$.
	
	\textbf{Case III:} There exist reactions $y_{j_1}X_{j_1} \to \tilde{y}'_{j_1\tilde{i}}X_{j_1}$ in $\tilde{\mathcal{M}}$, $K^{(y)}<0$ and $Y_{j_1}^{(y)}>0$. This means $\sum_{\tilde{y}'_{j_1\tilde{i}}>0}\bar{k}_{\tilde{i}}^{(y)}l_{j_1}\tilde{y}'_{j_1\tilde{i}}>Y^{(y)}_{j_1}$. Modify $\delta_{j_1}^{(y,\tau)}$ as
$\delta_{j_1}^{(y,\tau)}=Y^{(y,\tau)}_{j_1}/\sum_{\tilde{y}'_{j_1i}>0}\bar{k}_{\tilde{i}}^{(y)}l_{j_1}\tilde{y}'_{j_1 \tilde{i}}$ which is different with $\delta_j^{(y,\tau)}$ for $j\neq j_1$ defined as before. Based on this definition, $\delta_{j_1}^{(y)}=\sum_{\tau\in\bm{\tau}}\delta_{j_1}^{(y,\tau)}$ is less than 1 for $j_1$ and $\delta_{j}^{(y)}=1$ for $j\neq j_1$. Further, we denote $\bar{k}_{\tilde{i}}^{(y,\tau)}=\delta_{j_1}^{(y,\tau)}\bar{k}^{(y)}_{\tilde{i}}$ and $\bar{k}_{\tilde{i}}^{(y,0)}=(1-\delta_{j_1}^{(y)})\bar{k}^{(y)}_{\tilde{i}}$ which are the corresponding reaction rate constants of $y_{j_1}X_{j_1} \to \tilde{y}'_{j_1\tilde{i}}X_{j_1}$ with different  time delays $\tau$ and zero respectively. Then $$\sum_{\tau\in\bm{\tau}}Y^{(y,\tau)}x^y(t-\tau)=\sum_{\tau\in\bm{\tau}}\sum_{\tilde{i}=1}^{\tilde{r}}\bar{k}_{\tilde{i}}^{(y,\tau)}x^y(t-\tau)L\tilde{y}'_{.\tilde{i}},$$ $f^{(y)}$ can be written as $f^{(y)}=f_9+f_{10}$ where 
\begin{equation*}\label{eq:c3}
\begin{split}\small
f_9=&\sum_{\tau\in\bm{\tau}}Y^{(y,\tau)}x^y(t-\tau)+\sum_{\tilde{y}'_{j_1{\tilde{i}}>0}}\bar{k}_{\tilde{i}}^{(y,0)}x^y(t)L\tilde{y}'_{.i}-\sum_{\tilde{i}=1}^{\tilde{r}}\bar{k}_{\tilde{i}}^{(y)}x^y(t)Ly\\
		\end{split}
\end{equation*}	
\begin{equation*}
\begin{split}
&f_{10}=-\sum_{\tilde{y}'_{j_1{\tilde{i}}>0}}\bar{k}_{\tilde{i}}^{(y,0)}x^y(t)L\tilde{y}'_{.i}+(\sum_{\tilde{i}=1}^{\tilde{r}}\bar{k}_{\tilde{i}}^{(y)}l_{j_1}-\sum_{i=1}^{r}k_i^{(y)})x^y(t)y.
		\end{split}
	\end{equation*}
	Note that 
 \begin{equation*}
 \begin{split}
 &-\sum_{\tilde{y}'_{j_1{\tilde{i}}>0}}\bar{k}_{\tilde{i}}^{(y,0)}L\tilde{y}'_{.i}=(\delta_{j_1}^{(y)}-1)\sum_{\tilde{y}'_{j_1\tilde{i}}>0}\bar{k}_{\tilde{i}}^{(y)}L\tilde{y}'_{.i}=Y^{(y)}_{j_1}-\sum_{\tilde{y}'_{j_1\tilde{i}}>0}\bar{k}_{\tilde{i}}^{(y)}L\tilde{y}'_{.i}=K^{(y)}y.
 \end{split}
 \end{equation*}
 Thus $f_{10}$ is exactly 0 from above equation and the form of $K^{(y)}$. $f_9$ 
	 is the dynamics of all reactions whose reactant complex is $y$ in some $\ell$cQuasi-DCBMAS denoted as $\mathcal{M}_{\text{ql}}$ relative to  $\mathcal{\tilde{M}}$.
  
 \textbf{Case IV:} Consider the case that  $y=y_{j_1}e_{j_1}$ for some $j_1$ in $\mathcal{RC}$ but not in $\mathcal{R\tilde{C}}$. \eqref{eq:42} can be reduced into $Y^{(y)}=\sum_{i=1}^rk_i^{(y)}y_{j_1}e_{j_1}=Y_{j_1}^{(y)}e_{j_1}=\sum_{\tau\in\bm{\tau}}Y_{j_1}^{(y,\tau)}e_{j_1}$. Then the corresponding $f^{(y)}$ can be written as 
\begin{equation*}
f^{(y)}=\sum_{\tau\in\bm{\tau}}f^{(y,\tau)}=\sum_{\tau\in\bm{\tau}}Y_{j_1}^{(y,\tau)}e_{j_1}x^y(t-\tau)-\sum_{i=1}^rk_i^{(y)}y_{j_1}x^y(t)e_{j_1}.
	\end{equation*}
It is obviously a combination of some supporting modules 1 defined in Section \ref{sec:3}.
 
 Combining $f^{(y)}$ in the above four cases,  $\dot{x}$ is the sum of the dynamics of some $\ell$cQuasi-DCBMAS of $\mathcal{\tilde{M}}$ as the supporting module 3 and some supporting module 1. Thus we can derive the Lyapunov stability of $\mathcal{M}_{\text{De3}}$ by using the linear combination of Lyapunov functional of Supporting modules defined in Section \ref{sec:3} based on the form of $f^{(y)}$ for all $y\in\mathcal{RC}$.
$\hfill\blacksquare$
\begin{example}\label{ex:4} p21-activated kinase 1(PAK-1) network denoted as $\mathcal{M}_{\mathrm{De}}=(\mathcal{M}, \bm{\tau})$ plays a pivotal role in various biological processes such as tumor cell proliferation, apoptosis, and invasion \cite{H2021}.  $\mathcal{M}$ is linearly conjugate to $\mathcal{\tilde{M}}$ where $\tilde{k}_1=k_1/4$, $\tilde{k}_2=k_2$, $\tilde{k}_3=k_3$, $\tilde{k}_4=k_4$, and  $L=\mathrm{diag}(1/2,1,1)$ 
\begin{figure}[h]
	\centering
\includegraphics[height=3.5cm,width=9cm]{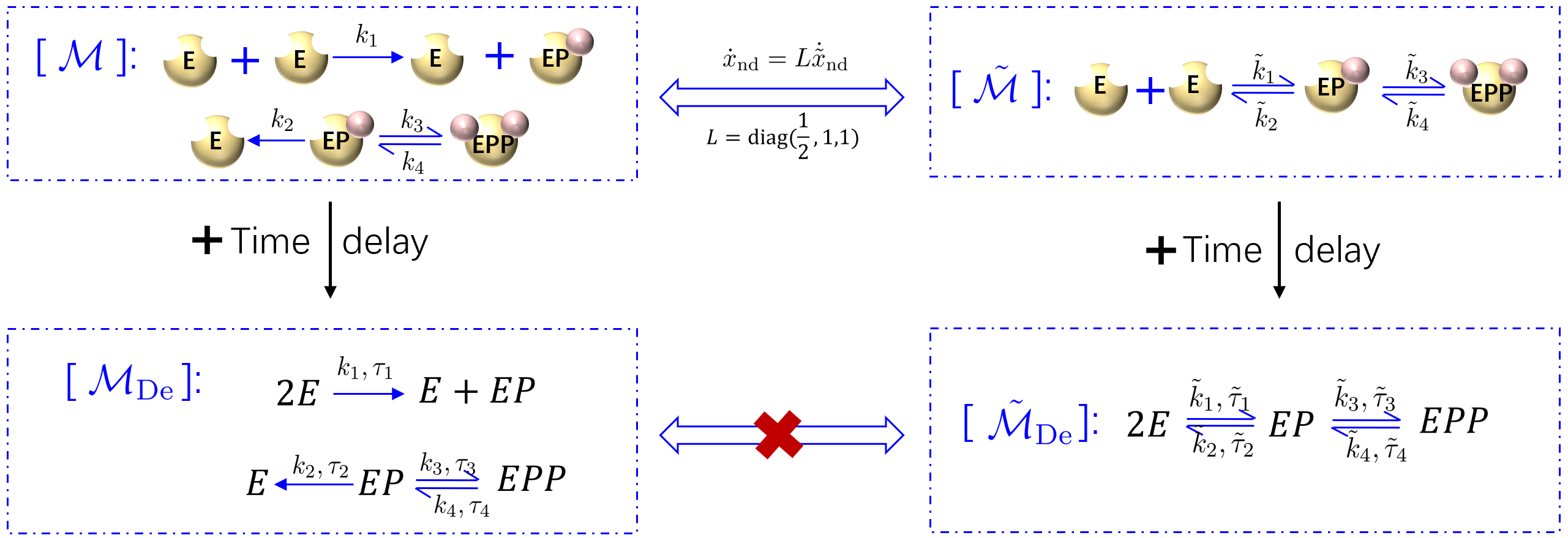}
	\caption{Delayed and non-delayed PAK-1 network and its corresponding $\mathcal{\tilde{M}}$.}
	\label{fig:PAK}
\end{figure}
The dynamics of $\mathcal{M}_{\mathrm{De}}$ is in the following form
\begin{equation*}
	\begin{split}
	\dot{x}&=k_1x_E^2(t-\tau_1) \left(\begin{matrix}
			1\\1\\0
		\end{matrix} \right)-k_1x_E^2(t)\left(\begin{matrix}
			2\\0\\0
		\end{matrix} \right)+k_2x_{EP}(t-\tau_2)\left(\begin{matrix}
			1\\0\\0
	\end{matrix}\right)-k_2x_{EP}(t)\left(\begin{matrix}
			0\\1\\0
		\end{matrix}\right)\\&+k_3x_{EP}(t-\tau_3)\left(\begin{matrix}
			0\\0\\1
		\end{matrix}\right)-k_3x_{EP}(t)\left(\begin{matrix}
		0\\1\\0
		\end{matrix}\right)+k_4x_{EPP}(t-\tau_4)\left(\begin{matrix}
		0\\1\\0
		\end{matrix}\right)-k_4x_{EPP}(t)\left(\begin{matrix}
		0\\0\\1
		\end{matrix}\right)
	\end{split}
\end{equation*}
which is not linear conjugate to $\mathcal{\tilde{M}}_{\mathrm{De}}$ for any $\tilde{\tau}_i$.
However, it is obviously a $\mathcal{M}_{\mathrm{De3}}$ and its dynamics can be written as 
\begin{equation*}
\dot{x}(t)=\dot{x}^{\mathrm{ql}}+k_1x^2_E(t-\tau_1)L\left(
	2, 0, 0
\right)^{\top}-k_1x^2_E(t)L\left(
	2, 0, 0
\right)^{\top}
\end{equation*}
from the description of Theorem \ref{thm:ss}. Then $\mathcal{M}_{\mathrm{De}}$ is Lyapunov stable.

\end{example}
This section outlines several methods to construct the Lyapunov functional for the general DeMAS $\mathcal{M}_{\text{De}}$. However, it's important to note that these methods primarily ensure the Lyapunov stability of the equilibrium. Local asymptotic stability, however, may be compromised due to the presence of degenerate equilibrium points within each compatibility class. Consequently, a more comprehensive analysis is required to determine the locally asymptotic stable of equilibrium points within SCC.
\section{The Preservation of Local Asymptotic Stability} \label{sec:5}
In this section, we demonstrate that the local asymptotic stability of D$\ell$cCBMASs cannot be guaranteed concerning SCC, primarily due to the existence of degenerate equilibriums. However, it is upheld concerning a newly defined invariant set. Before discussion, we revisit the relationship between $\mathscr{S}$ and $\tilde{\mathscr{S}}$ of $\mathcal{M}_{\text{De}}$ and $\mathcal{\tilde{M}}$, respectively.

 \begin{theorem}\cite{Feinberg1977}\label{thm:fh}
 Consider a MAS $\mathcal{M}$ with dynamics $\dot{x}_{\mathrm{nd}}=f(x)$ where the subspace spanned by the images of $f$ is referred to as the kinetics subspace $\mathscr{S}_k$. Then there exists $\mathscr{S}_k\subseteq \mathscr{S}$ and $\mathscr{S}_k$ is coincidence to $\mathscr{S}$ when $\mathcal{M}$ is a CBMAS $\tilde{\mathcal{M}}$. 
 \end{theorem}
 \begin{remark}\label{rem:5}
Suppose $\mathcal{M}_{\text{De}}=(\mathcal{M},\bm{\tau})$ is a D$\ell$cCBMAS.  Here, $\mathscr{S}_k$ and $\tilde{\mathscr{S}}_k$ are the kinetics subspace of $\mathcal{M}$ and the corresponding CBMAS $\mathcal{\tilde{M}}$ respectively, then there exists
$\mathscr{S}_k = L\mathscr{\tilde{S}}_k$. Additionally, combining $\mathscr{S}_k \subseteq \mathscr{S}$ and $\tilde{\mathscr{S}} = \tilde{\mathscr{S}}_k$ from Theorem \ref{thm:fh}, there exists $L\tilde{\mathscr{S}} \subseteq \mathscr{S}$. Further, dynamic equivalence implies $\tilde{\mathscr{S}} \subseteq \mathscr{S}$.
 \end{remark}
\begin{lemma}[Degenerate equilibriums in SCCs]\label{lem:de}
Each SCC of $\mathcal{M}_{\mathrm{De}}$ in Section \ref{sec:4} contains degenerate equilibriums if $\mathrm{dim}{\mathscr{S}}>\mathrm{dim}{\mathscr{\tilde{S}}}$.
 \end{lemma}
 \textit{\textbf{Proof:}}
 Denoting the dimensions of $\mathscr{S}$ and $\mathscr{\tilde{S}}$ as $s$ and $\tilde{s}$ respectively, we define $\textbf{b}\triangleq \{b_1,\cdots, b_{n-s}\}$ as a set of basis vectors for $\mathscr{S}^{\bot}$. If $s=n$, then $\mathcal{M}_{\text{De}}$ has only one SCC which is $\mathscr{C}_+$. Further $\tilde{s}<s=n$ indicates the existence of infinitely many positive equilibria in the phase space $\mathscr{C}_+$ of $\mathcal{M}_{\text{De}}$, as per the dynamics equivalence. Furthermore, the set of all equilibria can be expressed as $\{\bar{x}\vert \text{diag}(x^*)e^v, v\in \mathscr{\tilde{S}}^{\bot}\}$, where $x^*$ is some positive equilibrium. Thus, positive equilibria of $\mathcal{M}_{\text{De}}$ are degenerate.
 	
 If $\mathcal{M}_{\text{De}}$ is introduced in Subsection \ref{subsec:4.1} and \ref{subsec:4.2},  we consider the case where $s < n$, indicating that there are infinitely many SCCs in both $\mathcal{M}_{\text{De}}$ and $\tilde{\mathcal{M}}$. Each SCC $\mathcal{D}_{\psi}$ of $\mathcal{M}_{\text{De}}$ can be written as:
$
\mathcal{D}_{\psi}=\{\theta\vert c_{b}(\theta)=c_{b}(\psi)=W_b,~\text{for~each}~b\in \textbf{b}\}.$
  Further, for the case $\mathcal{RC}$ is equal to $\mathcal{R\tilde{C}}$, we introduce $\mathcal{\tilde{M}}_{\text{De}}=(\tilde{\mathcal{M}},\bm{\tilde{\tau}})$ where the time delay of each reaction with reactant complex $y$ denoted as $\tilde{\tau}^{(y)}$ is given by $\tilde{\tau}^{(y)}=\frac{\sum_{\tau\in\bm{\tau}}\sum_{i=1}^rk_{i}^{(y,\tau)}\tau}{\sum_{\tilde{i}=1}^{\tilde{r}}\tilde{k}_{\tilde{i}}^{(y)}}$. Thus, its SCC is denoted as $\mathcal{\tilde{D}}_{\psi}=\{\theta\vert \tilde{c}_{b}(\theta)=\tilde{c}_{b}(\psi)=W_b, b\in\bf{\tilde{b}}\}$.
  And for each constant function $\psi\in\bar{\mathscr{C}}_+$ with $\psi(s)=\psi, s\in[-\tau_{\text{max}},0]$ and each $b\in\textbf{b}$, there exists
  \begin{equation*}\label{eq:a}
  	\begin{split}
  &	c_{b}(\psi)=b^{\top}\left[\psi+\sum_{y\in\mathcal{RC}}\sum_{\tau\in\bm{\tau}}\sum_{i=1}^rk_i^{(y,\tau)}\psi^{y}\tau y\right]=b^{\top}\left[\psi+\sum_{y\in\mathcal{RC}}\sum_{\tilde{i}=1}^{\tilde{r}}\tilde{k}_{\tilde{i}}^{(y)}\psi^{y}\tilde{\tau}^{(y)} y\right]=\tilde{c}_{b}(\psi).
  	\end{split}
  \end{equation*}
  This reveals that all the constant functions in $\mathcal{D}_\psi$ are coincidence with those in $\bm{D}=\{\theta\vert \tilde{c}_b(\theta)=\tilde{c}_b(\psi)=W_b, b\in \textbf{b}\}$. 
$\mathscr{\tilde{S}}\subset \mathscr{S}$ ensures that $\textbf{b}$ is a subset of $\tilde{\bf{b}}$. Consequently, $\bm{D}$ is the union of infinite SCCs $\mathcal{\tilde{D}}_\psi$ and each such SCC in $\mathcal{\tilde{M}}_{\text{De}}$ contains a positive equilibrium. All these positive equilibriums are located within $\mathcal{D}_\psi$ of $\mathcal{M}_{\text{De}}$ and degenerate. If there exists $y\in \mathcal{RC}$ but not in $\mathcal{R\tilde{C}}$, we can add reactions $y\xrightarrow{\tilde{k}^{(y)},\tilde{\tau}^{(y)}} y$ into $\mathcal{\tilde{M}}_{\text{De}}$ where $\tilde{k}^{(y)}=\sum_{i=1}^rk_i^{(y)}$ and $\tilde{\tau}^{(y)}$ shares the same definition as above. Then we can conlcude the same conclusion. 

  Now consider $\mathcal{M}_{\text{De3}}$ relative to $\mathcal{\tilde{M}}$. We can introduce another DCBMAS $\tilde{\mathcal{M}}_{\text{De}}=(\tilde{\mathcal{M}},\tilde{\bm{\tau}})$ in which the time delay for each reaction with reactant complex $y$, denoted as $\tilde{\tau}^{(y)}$, is defined as $\frac{\sum_{i=1}^rk_i^{(y)}\sum_{j=1}^n\mathbbold{1}_{y_{ji}\neq 0}l^{-1}_{j}L^y\tau_i}{\sum_{\tilde{i}=1}^{\tilde{r}}\tilde{k}_{\tilde{i}}^{(y)}}$.
  Then for each $b\in\textbf{b}$, there exists
  \begin{equation}
  	\begin{split}
  		&c_b(\psi)=b^{\top}L(L^{-1}\psi+\sum_{i=1}^rk_i\psi^{y_{.i}}\tau_iL^{-1}y_{.i})\\
  		&=b^{\top}L\left[L^{-1}\psi+\sum_{y\in\mathcal{RC}}\sum_{i=1}^r(k_i^{(y)}L^{y}\sum_{j=1}^n\mathbbold{1}_{y_{ji}\neq0}l^{-1}_{j}\tau_i)(L^{-1}\psi)^{y}y\right]\\
  		&=b^{\top}L\left[L^{-1}\psi+\sum_{y\in\mathcal{RC}}\sum_{\tilde{i}=1}^{\tilde{r}}\tilde{k}_{\tilde{i}}^{(y)}\tilde{\tau}^{(y)}(L^{-1}\psi)^{y}y\right]\\
  		&=b^{\top}L\left[L^{-1}\psi+\sum_{\tilde{i}=1}^{\tilde{r}}\tilde{k}_{\tilde{i}}\tilde{\tau}_i(L^{-1}\psi)^{\tilde{y}_{.i}}\tilde{y}_{.i}\right]=\tilde{c}_{Lb}(L^{-1}\psi).
  	\end{split}
  \end{equation}
  where $\tilde{c}_{\tilde{b}}(L^{-1}\psi)$ represents one of functions in SCC $\mathcal{\tilde{D}}_{L^{-1}\psi}=\{\theta\vert \tilde{c}_{\tilde{b}}(\theta)=\tilde{c}_{\tilde{b}}(L^{-1}\psi), \tilde{b}\in \mathscr{\tilde{S}}^{\bot}\}$. The fact that $L\mathscr{\tilde{S}}$ is a subset of $\mathscr{S}$, as remarked in Remark \ref{rem:5}, implies that $\mathscr{S}^{\bot}$ is a subset of 
  $(L\mathscr{\tilde{S}})^{\bot}=\{L^{-1}\tilde{b}, \tilde{b}\in\mathscr{\tilde{S}}^{\bot}\}$. Consequently, $Lb$ corresponds to some $\tilde{b}$ in $\mathscr{\tilde{S}}^{\bot}$. Since $L\mathscr{S}^{\bot}$ is a subset of $\mathscr{\tilde{S}}$, we can derive an infinite number of degenerate equilibria in each $\mathcal{D}_{\psi}$, similar to the case discussed earlier.   
 	$\hfill\blacksquare$
 \begin{example}\label{ex:5}
Consider $\mathcal{M}_{\mathrm{De1}}$ relative to $\mathcal{\tilde{M}}$  and $\tilde{k}_1=\tilde{k}_2=k_1=k_2=k_3=k$.
 		\begin{figure}[H]
 			\centering
 			\includegraphics[height=1.1cm,width=8.9cm]{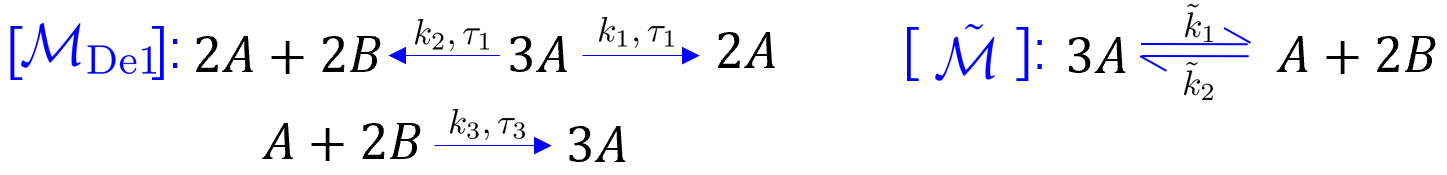}
 			\caption{A simple $\mathcal{M}_{\mathrm{De1}}$ with the existence of degenerate equilibriums.}
 			\label{fig:51}
 		\end{figure}
 $\mathscr{S}$ of $\mathcal{M}_{\mathrm{De1}}$ is 2d, which  equals the number of species. Thus, $\mathscr{S}^{\bot}$ is empty, indicating that there exists only one SCC — $\bar{\mathscr{C}}_+$. Each positive equilibrium $\bar{x}$ satisfying $\bar{x}_1=\bar{x}_2$  is also an equilibrium of $\mathcal{M}_{\mathrm{De1}}$. Thus, $\mathcal{M}_{\mathrm{De1}}$ exhibits degenerate equilibriums.
 
 The same situation can be found in the following $\mathcal{M}_{\mathrm{De3}}$ which is slightly different from the PAK-1 network as shown below:		
 	\begin{figure}[H]
 	\centering
\includegraphics[height=7cm,width=10cm]{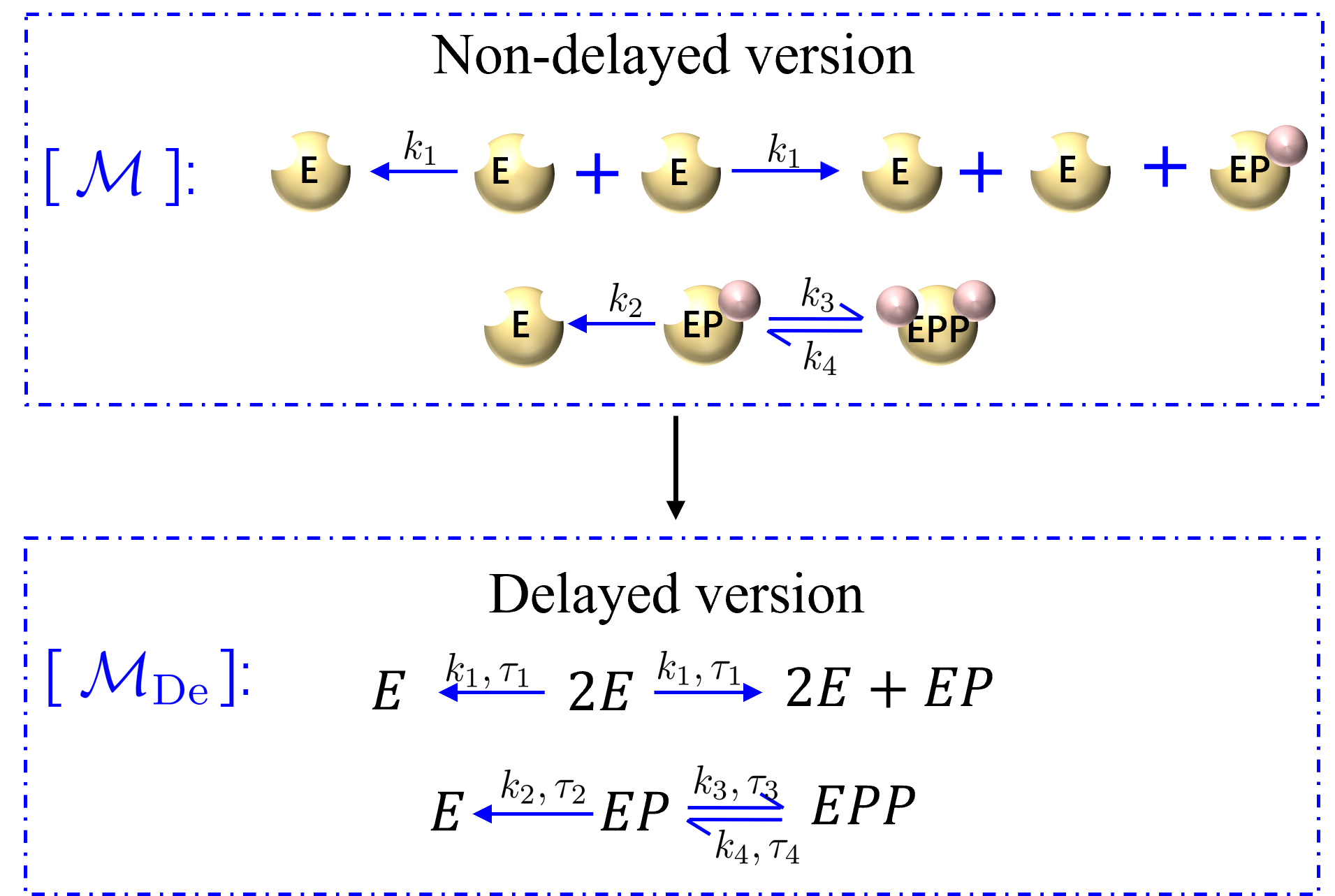}
 	\caption{A slight change of PAK-1 system which possesses degenerate equilibriums.}
 	\label{fig:scPAK}
 	 \end{figure}
	\noindent  where $k_1=k_2=k_3=k_4=\tilde{k}_2=\tilde{k}_3=\tilde{k}_4=k$ and $\tilde{k}_1=k/4$.  	$\mathcal{M}$ is linear conjugate to $\mathcal{\tilde{M}}$ defined in Example \ref{ex:4} with $L=\text{diag}(1/2,1,1)$. However, the SCC of $\mathcal{M}_{\mathrm{De3}}$ is $\bar{\mathscr{C}}_+$ with degenerate equilibriums  $\bar{x}_E^2=\bar{x}_{EP}=\bar{x}_{EPP}$ within it.
 	\end{example}
 
 \begin{remark}
 	The SCC outlined in Definition \ref{def:SCC} is established on the first integral of $\dot{x}$. These SCCs collectively form a partition of $\bar{\mathscr{C}}_+$ into infinitely equivalent classes, with each class serving as an invariant set of the trajectory. However, the first integral of $\dot{x}$ is not unique, resulting in various partitions. By engaging in re-decomposition, we can delve into the local asymptotic stability concerning other equivalent classes.
 \end{remark}
 \begin{lemma}\label{lem:nSCC}
 For each $\mathcal{M}_{\mathrm{De}}$ in Subsection \ref{subsec:4.1}, \ref{subsec:4.2}, and $\psi\in\bar{\mathscr{C}}_+$, the following $\mathcal{D}^n_\psi$  
 \begin{equation}
\mathcal{D}^n_\psi\triangleq\{\theta\vert c_{b}(\theta)=c_b(\psi)~\mathrm{for~each}~b\in \mathscr{\tilde{S}^{\bot}}\}
 \end{equation}
 where  $\mathscr{\tilde{S}}$ is the stoichiometric subspace of the $\mathcal{\tilde{M}}$ is a closed, positively invariant set under the dynamics of $\mathcal{M}_{\mathrm{De}}$.
 \end{lemma}
 \textit{\textbf{Proof:}}
 	The closeness of $\mathcal{D}^n_\psi$ as defined above is evident. Therefore, to conclude this lemma, it suffices to prove that each $c_b$ serves as the first integral of the dynamics of $\mathcal{M}_{\text{De}}$, indicating that the Lie derivative of each $c_b$ along every trajectory is zero. Since $\mathscr{S}^{\bot}$ is a subset of $\tilde{\mathscr{S}}^{\bot}$, we only need to verify the Lie derivative of $c_b$ where $b$ belongs to $\mathscr{S}^{\bot}$ but not in $\mathscr{\tilde{S}}^{\bot}$. For each $b$ of this kind, there exists
 \begin{equation*}
 	\begin{split}
L_{\dot{x}}c_b&=\dot{c}_b(x(t))
 =b^{\top}(\sum_{y\in\mathcal{RC}}\sum_{\tau\in\bf{\tau}}\sum_{i=1}^rk_i^{(y,\tau)}v_{.i}^{(y)}x^y(t-\tau)).
 	\end{split}
 \end{equation*}
  $\mathcal{M}_{\text{De}}$ provided in Subsections \ref{subsec:4.1} and \ref{subsec:4.2} indicates that  $\sum_{i=1}^{r}k_i^{(y,\tau)}v_{.i}^{(y)}$ for each $y$ can be expressed as a non-negative combination of some reaction vectors in $\tilde{\mathcal{M}}$, namely, $\sum_{\tilde{i}=1}^{\tilde{r}}e_{\tilde{i}}^{(y,\tau)}\tilde{v}^{(y)}_{.\tilde{i}}$. Since $b$ is a vector in $\mathscr{\tilde{S}}^{\bot}$, $L_{\dot{x}}c_b=0$ for each $b$ can be concluded.
 $\hfill\blacksquare$
 
 \begin{theorem}\label{thm:eue}
 Consider a	$\mathcal{M}_{\mathrm{De}}$ introduced in subsections \ref{subsec:4.1} and \ref{subsec:4.2}. Then the existence and uniqueness of the positive equilibrium within each $\mathcal{D}^n_{\psi}$ of $\mathcal{M}_{\mathrm{De}}$ holds. 
 \end{theorem}
 \textit{\textbf{Proof:}}
 	 Each $\mathcal{D}^n_{\psi}$ shares the same constant functions $\psi\in\bar{\mathscr{C}}_+$ with those in SCC $\mathcal{\tilde{D}}_{\psi}$ of DeCBMAS $\mathcal{\tilde{M}}_{\text{De}}$ defined in Lemma \ref{lem:de}. Thus, the existence and uniqueness of the equilibrium within each $\mathcal{D}^n_{\psi}$ can be derived from that in each SCC of $\mathcal{\tilde{M}}_{\text{De}}$.
 $\hfill\blacksquare$
 
  Now we consider the scenario where $\mathcal{M}_{\text{De3}}$ is described in Subsection \ref{subsec:4.3}.
 \begin{lemma}\label{lem:nis}
 For each $\mathcal{M}_{\mathrm{De3}}$, the $\mathcal{D}^n_\psi$ defined as 
\begin{equation}\label{eq:55}
\mathcal{D}^n_\psi\triangleq 
 \{\theta\vert c_b^u(\theta)=c^u_b(\psi)~\mathrm{for~each}~b\in (L\tilde{\mathscr{S}})^{\bot} \}
 \end{equation} 
is an invariant set of each trajectory where $$c_b^u(\psi)=b^{\top}[\psi(0)+\sum_{i=1}^r\int^0_{-\tau_i}k_i\sum_{j=1}^n\delta_{ji}x^{y_{.i}}(s)dsy_{.i}]\triangleq b^{\top}g^n(\psi)$$ and $\delta_{ji}=\frac{y'_{ji}}{\overline{l\tilde{y}'}_j}$ with $\overline{l\tilde{y}'}_j=\sum_{\tilde{y}'_{j\tilde{i}}\neq 0}\bar{k}_{\tilde{i}}^{(y)}l_j\tilde{y}'_{j\tilde{i}}/l_{j_1}\sum_{\tilde{y}'_{j\tilde{i}}\neq 0}\bar{k}_{\tilde{i}}^{(y)}$ for $j\neq j_1$ and $\overline{l\tilde{y}'}_{j_1}=y_{j_1}$.
\end{lemma}
 \textit{\textbf{Proof:}}
The Lie derivative of $c_b^u$ along each trajectory $x(t)$ is shown as below
\begin{equation}\label{eq:5.6}
	\begin{split}
		L_{\dot{x}}c_b^u&=\dot{c}_b(x(t))=b^{\top}\sum_{y\in\mathcal{C}}\sum_{i=1}^r(k_i^{(y)}x^{y}(t-\tau_i)y'_{.i}-k_i^{(y)}x^{y}(t)y)\\
		&+b^{\top}\sum_{y\in\mathcal{C}}\sum_{i=1}^rk_i^{(y)}\sum_{j=1}^n\delta_{ji}(x^y(t)-x^y(t-\tau_i))y.\\
	\end{split}
\end{equation}

 Recall that $Z^{(y)}=\sum_{i=1}^rk_i^{(y)}(y'_{.i}-y)=\sum_{\tilde{i}=1}^{\tilde{r}}\bar{k}_i^{(y)}L(\tilde{y}'_{.\tilde{i}}-y)$. 
  Further denoting  $y=y_{j_1}e_{j_1}$, $Z_j^{(y)}$ satisfies 
 \begin{equation*}\label{eq:5.7}
 	\begin{split}
Z^{(y)}_j=\sum_{i=1}^rk_i^{(y)}y'_{ji}~for~j\neq j_1; Z^{(y)}_{j_1}=\sum_{i=1}^rk_i^{(y)}(y'_{j_1i}-y_{j_1}).
 	\end{split}
 \end{equation*}
 Indeed, if there exists a reaction $y_{j_1}X_{j_1}\to \tilde{y}'_{j_1\tilde{i}}X_{j_1}$ in $\mathcal{\tilde{M}}$, then the vector $\tilde{v}^{(y)}_{.\tilde{i}}=(\tilde{y}'_{j_1\tilde{i}}-y_{j_1})e_{j_1}$ belongs to $\mathscr{\tilde{S}}$. Consequently, if a reaction $\tilde{R}_{\tilde{i}_1}$ with reactant complex $y$ satisfies $\tilde{y}'_{.\tilde{i}_1}=\tilde{y}'_{j\tilde{i}_1}e_j$, then each $e_{j}$ is also included in $\tilde{\mathscr{S}}$.
 From the description provided, it's evident that $g^n(\psi)$ can be expressed as a linear combination of $\{e_j \vert Z^{(y)}_j\neq 0\}$, which forms a subset of $\tilde{\mathscr{S}}$. Besides, in this case, $L\tilde{\mathscr{S}}$ is coincidence to $\tilde{\mathscr{S}}$ inducing $L_{\dot{x}}c_b^u=0$.
 
 If each reaction with reactant complex $y$ of $\mathcal{\tilde{M}}$ is in the form $y_{j_1}X_{j_1}\to \tilde{y}'_{j\tilde{i}}X_j$ where $j\neq j_1$, $Z^{(y)}$ satisfies that
 \begin{equation}\label{eq:5.8}
 \begin{split}
&Z^{(y)}_j=\sum_{\tilde{y}'_{j\tilde{i}}\neq 0}\bar{k}_{\tilde{i}}^{(y)}l_j\tilde{y}'_{j\tilde{i}}=\sum_{\tilde{y}'_{j\tilde{i}}\neq 0}\bar{k}_{\tilde{i}}^{(y)}l_{j_1}\overline{l\tilde{y}'}_j~\text{for}~j\neq j_1; \\&	Z^{(y)}_{j_1}=-\sum_{\tilde{i}=1}^{\tilde{r}}\bar{k}_{\tilde{i}}^{(y)}l_{j_1}y_{j_1}.
\end{split}
 \end{equation}
 Thus there exists $Z^{(y)}_j=\sum_{\tilde{y}'_{j\tilde{i}}\neq 0}\bar{k}_{\tilde{i}}^{(y)}l_{j_1}\overline{l\tilde{y}'_{j}}=\sum_{i=1}^rk_i^{(y)}\delta_{ji}\overline{l\tilde{y}'_j}$ which induces 
 \begin{equation}\label{eq:yy}
 	 \sum_{\tilde{y}'_{j\tilde{i}}\neq 0}\bar{k}_{\tilde{i}}^{(y)}l_{j_1}=\sum_{i=1}^rk_i^{(y)}\delta_{ji},\text{for~each~}j\neq j_1.
 \end{equation}


 Now we return to \eqref{eq:5.6}, for each reaction $y\to y'_{.i}$, $k_i^{(y)}y'_{.i}$ can be written as $k_i^{(y)}\sum_{j=1}^ny'_{ji}e_{j}=k_i^{(y)}\sum_{j=1}^n\delta_{ji}\overline{l\tilde{y}'}_je_{j}$. Therefore, $L_{\dot{x}}c^u_b$ is in the form of
\begin{multline*}
L_{\dot{x}}c_b^u=b^{\top}\sum_{y\in\mathcal{C}}\sum_{i=1}^r\biggl[k_i^{(y)}\sum_{j=1}^n\delta_{ji}(\overline{l\tilde{y}'}_je_j-y)x^y(t-\tau_i)+(k_i^{(y)}\sum_{j=1}^n\delta_{ji}-k_i^{(y)})x^y(t)y\biggr].
 \end{multline*}
 Utilizing \eqref{eq:yy}, the second part in above equation can further written as 
 \begin{equation*}
 \begin{split}
&\sum_{i=1}^r(k_i^{(y)}\sum_{j=1}^n\delta_{ji}-k_i^{(y)})x^y(t)y=\sum_{i=1}^r(k_i^{(y)}\sum_{j\neq j_1}\delta_{ji}+k_i^{(y)}\sum_{j=j_1}\delta_{j_1i}-k_i^{(y)})x^y(t)y
\end{split}
 \end{equation*}
 \begin{equation*}
 \begin{split}
		&=(\sum_{j=1}^n\sum_{\tilde{y}'_{j\tilde{i}}\neq 0}\bar{k}_{\tilde{i}}^{(y)}Ly+\sum_{i=1}^rk_i^{(y)}y'_{j_1i}e_{j_1}-\sum_{i=1}^rk_i^{(y)}y)x(t)^y
 		\\&=(\sum_{\tilde{i}=1}^{\tilde{r}}\bar{k}_{\tilde{i}}^{(y)}Ly+\sum_{i=1}^rk_i^{(y)}y'_{j_1i}e_{j_1}-\sum_{i=1}^rk_i^{(y)}y)x(t)^y\\&=(Z^{(y)}_{j_1}-Z^{(y)}_{j_1})e_{j_1}x(t)^y=0
 	\end{split}
 \end{equation*}
 
Indeed, each $\overline{l\tilde{y}'_{j}}e_j-y$ can be expressed as $$\overline{l\tilde{y}'_{j}}e_j-y=\frac{\sum_{\tilde{y}'_{j\tilde{i}}\neq 0}\bar{k}_i^{(y)}(l_j\tilde{y}'_{.\tilde{i}}-l_{j_1}y)}{\sum_{\tilde{y}'_{j\tilde{i}}\neq 0}\bar{k}_i^{(y)}l_{j_1}}$$ where each $ l_j\tilde{y}'_{.\tilde{i}}-l_{j_1}y=L(\tilde{y}'_{.\tilde{i}}-y)$ is in $L\mathscr{\tilde{S}}$, thus the first part of $L_{\dot{x}}c_b^u$ is in $L\mathscr{\tilde{S}}$. Therefore, we can conclude that $L_{\dot{x}}c_b^u$ is zero from $b\in (L\mathscr{\tilde{S}})^{\bot}$. 

In the case where $y$ is in $\mathcal{RC}$ of $\mathcal{M}_{\text{De3}}$ but not in $\mathcal{R\tilde{C}}$ of $\tilde{\mathcal{M}}$, given that $Z^{(y)}=0$, we can deduce that each product complex $y'_{.i}$, with reactant complex $y=y_{j_1}e_{j_1}$, is also of the form $y'_{.i}=y'_{j_1i}e_{j_1}$. By defining $\overline{l\tilde{y}'}_{j_1}$ and $\delta_{j_1i}$ as specified in the conditions of this lemma, the derivation of $L_{\dot{x}}c_b^u=0$ can proceed analogously.
$\hfill\blacksquare$

The following Lemma aids proving  the existence and uniqueness of equilibrium.
  \begin{lemma}\cite{G2019}
 	Considering a linear subspace $S$ of $\mathbb{R}^n$ and $\eta_1, \eta_2\in\mathbb{R}_+^n$, $\beta_i\in\mathbb{R}_+$ and  $y_{.i}\in \mathbb{R}^n_{\geq 0}$ for each $1\leq i\leq r$. Then there must exist some $v\in S^{\bot}$ such that 
 	\begin{equation}\label{eq:5.12}
 		\text{diag}{(\eta_1)}e^v+\sum_{i=1}^r\beta_i\left[\text{diag}{(\eta_1)}e^v\right]^{y_{.i}}y_{.i}-\eta_2\in S.
 	\end{equation}
 \end{lemma}
 \begin{theorem}\label{thm:eu}
  For arbitrary $\mathcal{M}_\mathrm{De3}$ defined in Section \ref{sec:4}, the positive equilibrium is existed and unique within each $\mathcal{D}^n_\psi$.
 	\end{theorem}
 	\textit{\textbf{Proof:}}
\textbf{Existence:} 
The linear conjugacy between $\mathcal{M}_{\text{De3}}$ and $\mathcal{\tilde{M}}$ implies the existence of a positive equilibrium $\bar{x}\in \mathbb{R}_{>0}^n$ in $\mathcal{M}_{\text{De3}}$. Moreover, each equilibrium $x^*$ of $\mathcal{M}_{\text{De3}}$ can be expressed as $x^*=\text{diag}{(\bar{x})}e^v$, where $v$ belongs to $\mathscr{\tilde{S}}^{\bot}$. For each $\mathcal{D}^n_{\psi}$ of $\mathcal{M}_{\text{De3}}$ containing $\psi\in\bar{\mathscr{C}}+$, by setting $\eta_1=L^{-1}\bar{x}, S=\mathscr{\tilde{S}}, \beta_i=k_i\sum_{j=1}^n\delta_{ji}/\sum_{j=1}^n\mathbbold{1}_{y_{ji}\neq0}l_{j},  ~\mathrm{and}~ \eta_2=L^{-1}g^n(\psi),$we can assert that there exists a $v$ such that \eqref{eq:5.12} holds. Consequently, there exists an equilibrium $x^*=\text{diag}{(\bar{x})}e^v$ such that $L^{-1}g^n(x^*)-L^{-1}g^n(\psi)\in \mathscr{\tilde{S}}$, thus confirming the existence of the positive equilibrium $x^*$ in $\mathcal{D}_{\psi}^n$.
	
\textbf{Uniqueness:} Assuming $x^*$ and $x^{**}$ are two equilibria in the same invariant set $\mathcal{D}^n_{\psi}$, they satisfy  $\ln{x^*}-\ln{x^{**}}\in \mathscr{S}^{\bot}$ and 
\begin{multline*}
	(\ln{x^*}-\ln{x^{**}})^{\top}L^{-1}\biggl[L(L^{-1}x^*-L^{-1}x^{**})+\\
	\sum_{y\in\mathcal{RC}}\sum_{i=1}^rk^{(y)}_i/l^{(y)}_{j_1}\sum_{j=1}^n\delta_{ji}\tau_i[(x^*)^{y}-(x^{**})^{y}]Ly\biggr]=0
\end{multline*}
Further combining $(\ln{x_1}-\ln{x_2})(x_1-x_2)=0$ iff $x_1=x_2$, we can obtain $x^*=x^{**}$. 
	$\hfill\blacksquare$

 \begin{theorem}\label{thm:las}
	Each $\mathcal{M}_{\mathrm{De}}$ described in Section \ref{sec:4} is locally asymptotically stable relative to invariant classes $\mathcal{D}^n_{\psi}$.  
\end{theorem}
\textit{\textbf{Proof:}}
	It can be derived from the results in Section \ref{sec:4} and the existence, uniqueness of the equilibrium relative to $\mathcal{D}_{\psi}^n$ shown in Theorems \ref{thm:eue} and \ref{thm:eu}. 
$\hfill\blacksquare$

\begin{example}
Now we explore the newly defined invariant sets of systems $\mathcal{M}_{\mathrm{De}}$ presented in Fig. \ref{fig:51} and \ref{fig:scPAK} respectively. 
	
	For the system depicted in Fig. \ref{fig:51}, each invariant set $\mathcal{D}^n_\psi$, where $\psi\in\mathscr{S}$, shares constant functions with the corresponding stoichiometric compatibility class $\tilde{\mathcal{D}}_\psi$ of the DCBMAS $\tilde{\mathcal{M}}$. This correspondence ensures that the existence and uniqueness of equilibria in $\tilde{\mathcal{M}}_{\mathrm{De}}$ are preserved in $\mathcal{M}_{\mathrm{De}}$. Moreover, the local asymptotic stability of $\mathcal{M}_{\mathrm{De}}$ can be inferred from Theorem \ref{thm:las}.

	Now, we further consider $\mathcal{M}_{\mathrm{De}}$ as depicted in Fig. \ref{fig:scPAK}. Its non-delayed version is linear conjugate to $\mathcal{\tilde{M}}$, with each complex being a multiple of some single species. Without loss of generality, we denote species $E$, $EP$, and $EPP$ as $X_j$, where $j=1,2,3$ respectively.
		The newly defined invariant set of $\mathcal{M}_{\text{De}}$, denoted as $\mathcal{D}^n_{\psi}$ for each $\psi\in\bar{\mathscr{C}}_{+}$, is shown in \eqref{eq:55}, where $\delta_{ji}$ inside $c^u(\psi)$ can be computed as shown in Lemma \ref{lem:nis}.  
	\begin{itemize}
		\item For the reactant complex $2E$,  $\overline{l\tilde{y}}_{1}=y_1=2$, $\overline{l\tilde{y}}_2=1/l_{1}=2$, $\overline{l\tilde{y}}_3=0$. Further the corresponding $\delta_{ji}$ are $\delta_{11}=1/2$, $\delta_{21}=0$, $\delta_{31}=0$ and $\delta_{12}=2/2=1$, $\delta_{22}=1/2$, $\delta_{32}=0$.
		\item For the reactant complex $EP$, $\overline{l\tilde{y}}_{1}=2*(1/2)/1=1$, $\overline{l\tilde{y}}_{2}=y_2=1$, $\overline{l\tilde{y}}_{3}=1/1=1$. Further, the corresponding $\delta_{ji}$ are 
	$\delta_{13}=1$, $\delta_{23}=0$, $\delta_{33}=0$ and $\delta_{14}=0$, $\delta_{24}=0$, $\delta_{34}=1$.
	\item For the reactant complex $EPP$, $\overline{l\tilde{y}}_{1}=0$, $\overline{l\tilde{y}}_{2}=1$, $\overline{l\tilde{y}}_{3}=1$. Further, the corresponding $\delta_{ji}$ are 
	$\delta_{15}=0$, $\delta_{25}=1$, $\delta_{35}=0$. 
	\end{itemize}

	Now we verify that $\mathcal{D}_\psi^n$ is  exactly an invariant set of $\mathcal{M}_{\mathrm{De}}$. We can verify that $g^n(x)$ is
\begin{equation*}
	\begin{split}
&g^n(x)=x(t)+\sum_{i=1}^r\int^0_{-\tau_i}k_i\sum_{j=1}^n\delta_{ji}x^{y_{.i}}(t+s)dsy_{.i}\\
	&=k_1x_1^2(t-\tau_1)\left(\begin{matrix}
			-1\\1\\0
		\end{matrix}\right)+k_2x_2(t-\tau_2)\left(\begin{matrix}
			1\\-1\\0
		\end{matrix}\right)+k_3x_2(t-\tau_3)\left(\begin{matrix}
			0\\-1\\1
		\end{matrix}\right)+k_4x_3(t-\tau_4)\left(\begin{matrix}
			0\\1\\-1
		\end{matrix}\right).
	\end{split}
\end{equation*}
 $\mathscr{\tilde{S}}$ is spanned by the reaction vectors $\{(-2,1,0), (0,-1, 1)\}$. It is easy to verify that $L_{\dot{x}}c^u_b$ equals zero, as each vector in the last line of the above equation belongs to $L\mathscr{\tilde{S}}$.

By setting $k_i=1$ and $\tau_i=1$, Figure \ref{fig:52} illustrates the degenerate equilibriums $\bar{x}_1^2=\bar{x}_2=\bar{x}_3$ in the unique SCC of $\mathcal{M}_{\mathrm{De}}$ which is exactly $\bar{\mathscr{C}}_+$. Additionally, two surfaces contain all constant functions of the corresponding newly defined invariant sets $\mathcal{D}^n_{\psi_1}$ (the bigger one) and $\mathcal{D}^n_{\psi_2}$ respectively. Each invariant set intersects with the equilibrium set exactly once, with the unique equilibriums in $\mathcal{D}_{\psi_1}$ and $\mathcal{D}_{\psi_2}$ being $\bar{x}_1=(1.62,2.6245,2.6245)^{\top}$ and $\bar{x}_2=(0.8,0.64,0.64)^{\top}$ respectively.
Curves in blue and red represent different trajectories originating from four distinct initial functions: $\theta_1=(0.1,0.9,11.2)^{\top}$ (in blue) and $\theta_2=(2.2,0.7,0.79)^{\top}$ (in red) on $\mathcal{D}^n_{\psi_1}$. $\theta_3=(0.1,0.4,2.61)^{\top}$ (in blue) and  $(1.1,0.2,0.01)^{\top}$ (in red) on $\mathcal{D}^n_{\psi_2}$. They all converge to the unique equilibrium in the same invariant set.
 \begin{figure}[H]
	\centering
\includegraphics[height=10cm,width=13.0cm]{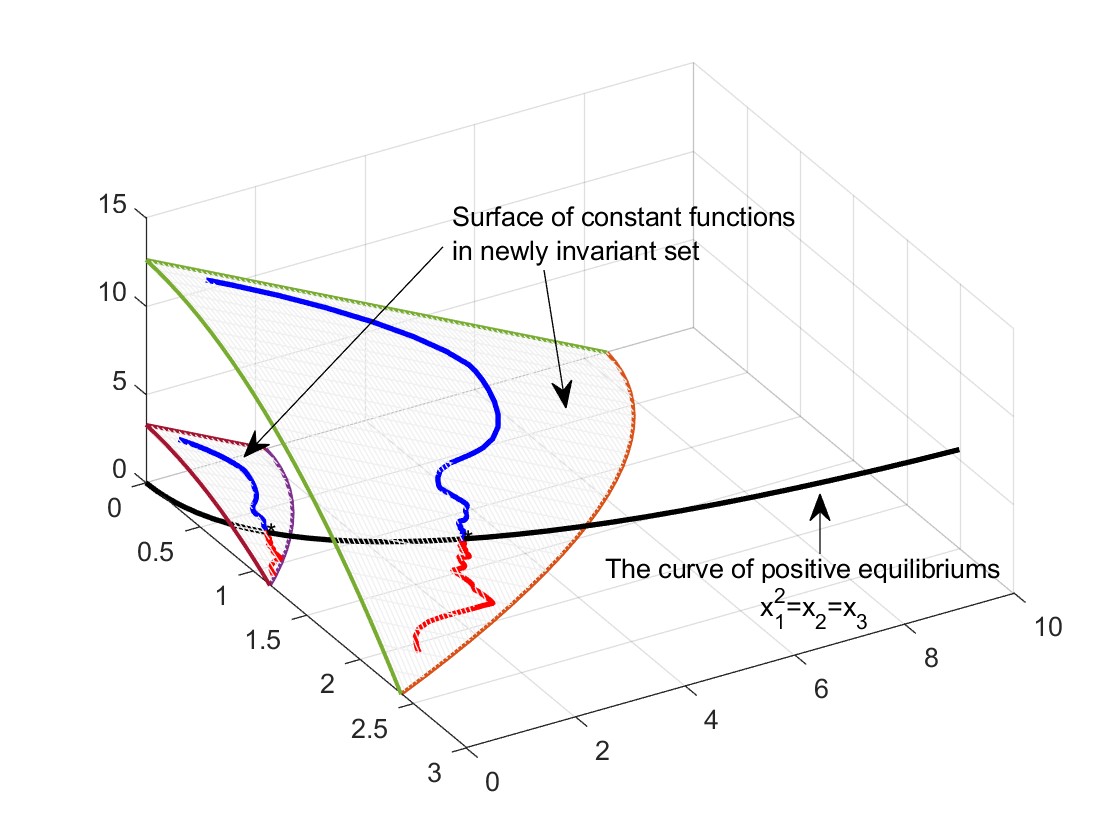}
	\caption{Two surfaces of constant functions of the newly defined invariant sets and the evolution of four trajectories on these sets.}
	\label{fig:52}
\end{figure}
	\end{example}
\begin{appendices}
\section{Lemma for the proof of Theorem \ref{thm:bz}}

\begin{lemma}\label{lem:ddx} Consider a $\mathcal{M}_{\mathrm{De}}=(\mathcal{M}, \bm{\tau})$ with $\mathcal{M}$ sharing the same dynamics as a CBMAS $\mathcal{\tilde{M}}$.  There exist sets $c^{(y,\tau,+)}$ and $c^{(y,\tau,-)}$ for each $y$ and $\tau$ that satisfy
\begin{equation}\label{eq:c}
 \begin{split}
&0<\sum_{i=1}^rk_i^{(y,\tau)}a^{(y,\tau)}_i\triangleq Z_+^{(y,\tau)}=\sum_{\tilde{a}^{(y)}_{\tilde{i}}>0}c^{(y,\tau,+)}_{\tilde{i}}\tilde{a}^{(y)}_{\tilde{i}}\\
 			&0>\sum_{i=1}^rk_i^{(y,\tau)}a^{(y,\tau)}_i\triangleq Z_-^{(y,\tau)}=\sum_{\tilde{a}^{(y)}_{\tilde{i}}<0}c^{(y,\tau,-)}_{\tilde{i}}\tilde{a}^{(y)}_{\tilde{i}}
 		\end{split}
 	\end{equation}
 	if both $\tilde{a}_{\tilde{i}_1}^{(y)}>0$ and $\tilde{a}_{\tilde{i}_2}^{(y)}<0$ exist for some $\tilde{i}_1$ and $\tilde{i}_2$.
\end{lemma}
\textit{\textbf{Proof:}}
	Dynamical equivalence ensures that $Z^{(y)}$ is the same for $\mathcal{M}$ and $\mathcal{\tilde{M}}$
	\begin{equation}
		Z^{(y)}=Z^{(y)}_+-Z^{(y)}_-=\tilde{Z}^{(y)}_+-\tilde{Z}^{(y)}_-
	\end{equation}
	where $Z^{(y)}_+, Z^{(y)}_-$ are the sum of all the absolute values of $Z^{(y,\tau)}>0$ and $Z^{(y,\tau)}<0$, respectively.
	And $\tilde{Z}^{(y)}_+$ and $\tilde{Z}^{(y)}_-$ can be defined as $\tilde{Z}^{(y)}_+=\sum_{\tilde{a}^{(y)}_{\tilde{i}}>0}\tilde{k}^{(y)}_{\tilde{i}}\tilde{a}^{(y)}_{\tilde{i}}$ and $\tilde{Z}^{(y)}_-=\sum_{\tilde{a}^{(y)}_{\tilde{i}}<0}\tilde{k}^{(y)}_{\tilde{i}}\vert\tilde{a}^{(y)}_{\tilde{i}}\vert$. Note that $\tilde{Z}^{(y)}_+$ and $\tilde{Z}^{(y)}_-$ may be different from $Z^{(y)}_+$ and $Z^{(y)}_-$. Further define
	$c^{(y,\tau,+)}_{\tilde{i}}$ and $c^{(y,\tau,-)}_{\tilde{i}}$ as follows
	\begin{equation}
		c^{(y,\tau,+)}_{\tilde{i}}\triangleq \delta_+^{(y,\tau)}\tilde{k}_{\tilde{i}};~ c^{(y,\tau,-)}_{\tilde{i}}\triangleq \delta_-^{(y,\tau)}\tilde{k}_{\tilde{i}}
	\end{equation}
	where $\delta_+^{(y,\tau)}=\frac{Z_+^{(y,\tau)}}{\tilde{Z}^{(y)}_+}$ and $\delta_-^{(y,\tau)}=\frac{\vert Z^{(y,\tau)}_-\vert}{\tilde{Z}^{(y)}_-}$. It is easy to verify that the coefficients we defined in the above equation satisfy the equation \eqref{eq:c}.
$\hfill\blacksquare$
\begin{remark}\label{rem:dd}
If each $Z^{(y,\tau)}$ and $Z^{(y)}$ have the same sign, without loss of generality, we suppose that they all larger than zero, namely, $Z^{(y)}_-=0$. Thus by setting $\delta^{(y,\tau)}=Z^{(y,\tau)}/Z^{(y)}$ and $\tilde{k}^{(y,\tau)}=\tilde{k}^{(y)}\delta^{(y,\tau)}$, we can obtain that 
$$
Z^{(y,\tau)}=\delta^{(y,\tau)}Z^{(y)}=\delta^{(y,\tau)}\tilde{Z}^{(y)}=\sum_{\tilde{i}=1}^{\tilde{r}}\delta^{(y,\tau)}\tilde{k}^{(y)}_{\tilde{i}}\tilde{a}^{(y)}_{\tilde{i}}=\sum_{\tilde{i}=1}^{\tilde{r}}\tilde{k}^{(y,\tau)}_{\tilde{i}}\tilde{a}^{(y)}_{\tilde{i}}.$$
\end{remark}
\section{Proof of Theorem \ref{thm:bz}}
\textit{\textbf{Proof:}}
(1) Case I: 
Firstly, we rewrite each $f^{(y,\tau)}$ as follows:  
$$
	f^{(y,\tau)}
=\sum_{i=1}^r\biggl[k_i^{(y,\tau)}x(t-\tau)^{y}y+k_i^{(y,\tau)}x(t-\tau)^{y}a_i^{(y,\tau)}w-k_i^{(y,\tau)}x(t)^{y}y\biggr]
$$
Combining Remark \ref{rem:dd} shown in Appendix A, we have 
$$
\dot{x}=\sum_{y\in\mathcal{RC}}\sum_{\tau\in\bm{\tau}}\biggl[\sum_{\tilde{i}=1}^{\tilde{r}}[\tilde{k}_{\tilde{i}}^{(y,\tau)}x(t-\tau)^{y}\tilde{y}'_{.i}-\tilde{k}_{\tilde{i}}^{(y,\tau)}x(t)^{y}y]+K^{(y,\tau)}(x(t-\tau)^{y}-x(t)^{y})y\biggr]
$$
\begin{equation*}
	\begin{split}
&=\dot{x}^{\text{q}}+\sum_{y\in\mathcal{RC}}\sum_{\tau\in\bm{\tau}}K^{(y,\tau)}(x(t-\tau)^{y}-x(t)^{y})y\triangleq \dot{x}^{\text{q}}+f_{3}
\end{split}
\end{equation*}
where $K^{(y,\tau)}\triangleq\sum_{i=1}^rk_i^{(y,\tau)}-\sum_{\tilde{i}=1}^{\tilde{r}}\tilde{k}^{(y,\tau)}_{\tilde{i}}$. 

The first part, denoted as $\dot{x}^{\text{q}}$, is exactly the dynamics of a Quasi-DeCBMAS relative to $\mathcal{\tilde{M}}$ and its Lyapunov functional is $V_{\text{q}}(\psi)$. Further from the definition of HF-DeMAS and the Remark \ref{rem:dd}, $K^{(y,\tau)}>0$ can be obtained. Supporting module 1 reveals the Lyapunov functional of $f_3$, denoted as $V_3(\psi)$. 
Thus, the Laypunov stability of each positive equilibrium $\bar{x}$ of $\mathcal{M}_{\text{De2}}$ can be derived by using the Lyapunov functional $V(\psi)=V_{\text{q}}(\psi)+V_3(\psi)-\sum_{j=1}^nh(\psi_j(0))$ and its Lie derivative along each trajectory $x(t)$ is zero which concludes the Lyapunov stability of $\mathcal{M}_{\text{De2}}$.
 
 (2) Case II: In this case, we can categorize all $f^{(y,\tau)}$ for each $y$ into two classes denoted as $f^{(y,\tau)}_-$ and $f^{(y,\tau)}_+$ based on $Z^{(y,\tau)}_-$ and $Z^{(y,\tau)}_+$ respectively. Consequently, the dynamics $\dot{x}$ of $\mathcal{M}_{\text{De2}}$
 can be expressed as:  
 $$\dot{x}=\sum_{y\in\mathcal{RC}}(\sum_{Z^{(y,\tau)}_+}f^{(y,\tau)}_++\sum_{Z^{(y,\tau)}_-}f^{(y,\tau)}_-).$$

By applying Lemma \ref{lem:ddx}, $f^{(y,\tau)}_+$ and $f^{(y,\tau)}_-$ can be expressed in the following form
\begin{equation*}
	\begin{split}
f^{(y,\tau)}_+&=\underbrace{\sum_{\tilde{a}^{(y)}_{\tilde{i}}>0} (c^{(y,\tau,+)}_{\tilde{i}}x(t-\tau)^y\tilde{y}'_{.\tilde{i}}-c^{(y,\tau,+)}_{\tilde{i}}x(t)^yy)}_{f^{(1)}_+}+\underbrace{K^{(y,\tau)}_+(x(t-\tau)^y-x(t)^y)y}_{f^{(2)}_+}\\
f^{(y,\tau)}_-&=\underbrace{\sum_{\tilde{a}^{(y)}_{\tilde{i}}<0} (c^{(y,\tau,-)}_{\tilde{i}}x(t-\tau)^y\tilde{y}'_{.\tilde{i}}-c^{(y,\tau,-)}_{\tilde{i}}x(t)^yy)}_{f^{(1)}_-}+\underbrace{K^{(y,\tau)}_-(x(t-\tau)^y-x(t)^y)y}_{f^{(2)}_-}
\end{split}
\end{equation*}
where $K_+^{(y,\tau)}=\sum_{i=1}^rk_i^{(y,\tau)}-\sum_{\tilde{a}^{(y)}_{\tilde{i}}>0}c^{(y,\tau,+)}_{\tilde{i}}$, $K_-^{(y,\tau)}=\sum_{i=1}^rk_i^{(y,\tau)}-\sum_{\tilde{a}^{(y)}_{\tilde{i}}<0}c^{(y,\tau,-)}_{\tilde{i}}$.

 $\dot{x}$ is in the following form
\begin{equation*}
\begin{split}
	\dot{x}=f^{(1)}+f^{(2)}=\sum_{y\in\mathcal{RC}}(f_+^{(1)}+f_-^{(1)})+\sum_{y\in \mathcal{RC}}(f_+^{(2)}+f_-^{(2)}).
\end{split}
\end{equation*}
 First consider $\tilde{Z}^{(y)}_+\geq Z^{(y)}_+$, further combining the fact $\tilde{Z}^{(y)}_+-Z^{(y)}_+=\tilde{Z}^{(y)}_--Z^{(y)}_-$, $f^{(1)}$ can be interpreted as
	 \begin{multline*}
	     f^{(1)}+\sum_{y\in\mathcal{RC}}\biggl[\sum_{\tilde{a}^{(y)}_{\tilde{i}}>0}(\tilde{k}_{\tilde{i}}^{(y)}-\sum_{\tau\in\bm{\tau}}c_{\tilde{i}}^{(y,\tau,+)})x^y(t)\tilde{v}^{(y)}_{.\tilde{i}}+\sum_{\tilde{a}^{(y)}_{\tilde{i}}<0}(\tilde{k}_{\tilde{i}}^{(y)}-\sum_{\tau\in\bm{\tau}}c_{\tilde{i}}^{(y,\tau,-)})x^y(t)\tilde{v}^{(y)}_{.\tilde{i}}\bigg]
      \end{multline*}
	 due to the part (a) is zero. 
In this context, $f^{(1)}$ also characterizes the dynamics of a Quasi-DCBMAS $\tilde{\mathcal{M}}_{\text{q}}$ relative to $\mathcal{\tilde{M}}$. Here, $\tilde{k}_{\tilde{i}} - \sum_{\tau \in \bm{\tau}} c_{\tilde{i}}^{(y,\tau,+)}$, denoted as $\tilde{k}^{(y,0)}_{\tilde{i}}$, represents the rate constant associated with the $\tilde{i}$-th reaction in $\mathcal{\tilde{M}}_{\text{q}}$ under zero time delay. The same interpretation applies to $\tilde{k}_{\tilde{i}} - \sum_{\tau \in \bm{\tau}} c_{\tilde{i}}^{(y,\tau,-)}$.
 The Lyapunov functional of $f^{(1)}$ can be derived from \textbf{supporting module 2} and is denoted as $V_{\text{q}}$. Furthermore, $$Z_{+}^{(y,\tau)}=\sum_{i=1}^rk_i^{(y,\tau)}a_i^{(y,\tau)}=\sum_{\tilde{\alpha}^{(y)}_{\tilde{i}}>0}c^{(y,\tau,+)}\tilde{a}_{\tilde{i}}^{(y,\tau)}.$$ As $K^{(y,\tau)}\geq 0$ from the Definition \ref{def:hf} of $\mathcal{M}_{\text{De2}}$,   the Lyapunov functional of $\dot{x}$ can be expressed similar to Case I.	

Then consider the scenario that there exists $y$ such that $Z^{(y)}_+$ and $Z^{(y)}_-$ are smaller than $\tilde{Z}^{(y)}_+$ and $\tilde{Z}^{(y)}_-$. The following part contributes to presenting that $f^{(1)}$ is exactly the dynamics of a Quasi-DeCBMAS $\tilde{\mathcal{M}}'_{\text{q}}=(\mathcal{\tilde{M}}', \bm{\tau}^{\text{q}})$ where the CBMAS $\mathcal{\tilde{M}}'$ shares the same network graph and positive equilibria with $\mathcal{\tilde{M}}$. 

To do this, we examine the network graphs $G=(V,E)$ and $G'=(V',E')$ of $\mathcal{\tilde{M}}$ and $\mathcal{\tilde{M}'}$ respectively.
Firstly, the sets of reactant complexes of $\mathcal{\tilde{M'}}$ and $\mathcal{\tilde{M}}$ are both $\mathcal{RC}$. Lemma \ref{lem:ddx} reveals the  positivity of both $c^{(y,\tau,+)}$ and $c^{(y,\tau,-)}$ which ensures the sets of product complexes are also the same for $\mathcal{\tilde{M}}$ and $\mathcal{\tilde{M}}'$  and each reaction $y \to \tilde{y}'_{.i}\in \tilde{\mathcal{M}}$ also exists in $\tilde{\mathcal{M}}'$. Now we consider the equilibrium of $\mathcal{\tilde{M}}'$. It is easy to verify that each constant function in $\mathscr{C}_+$ is an equilibrium of $f^{(2)}$. 
 Therefore, $\dot{\bar{x}}=f^{(1)}(x)+f^{(2)}(x)$ is equal to $f^{(1)}(x)$ for each constant function $x\in \mathscr{C}_+$. Thus we can conclude that the equilibriums of $\mathcal{\tilde{M}}'$ are coincidence to those of $\mathcal{\tilde{M}}$. And the Lyapunov functional of the Quasi-DeCBMAS $\tilde{\mathcal{M}}_{\text{q}'}$ can be denoted as $V_{\text{q}'}$.
 
 Now we further consider $f^{(2)}(x)$. 
 For each $y$ and $\tau$, there exists \begin{equation*}\begin{split}
&Z^{(y,\tau)}_+w^{(y)}=\sum_{i=1}^{r}k_i^{(y,\tau)}\bar{v}^{(y,\tau)}=\sum_{\tilde{a}_{\tilde{i}}^{(y)}>0}c^{(y,\tau,+)}_{\tilde{i}}\tilde{v}_{.\tilde{i}}^{(y,\tau)}=\sum_{\tilde{a}_{\tilde{i}}^{(y)}>0}c^{(y,\tau,+)}_{\tilde{i}}\bar{\tilde{v}}^{(y,\tau)}
\end{split}
\end{equation*}where $\bar{v}^{(y,\tau)}$ and $\bar{\tilde{v}}^{(y,\tau)}$ are the weighted averages of $v^{(y,\tau)}_{i}$ and $\tilde{v}^{(y,\tau)}_{\tilde{i}}$ relative to $k_i^{(y,\tau)}$ and $c^{(y,\tau,+)}_{\tilde{i}}$ respectively. Thus, $\Vert \bar{v}\Vert_1\leq\Vert \bar{\tilde{v}}\Vert_1$ can be obtained due to the Definition \ref{def:hf}, which further indicates $K_+^{(y,\tau)}\geq 0$.
 The same holds for $K^{(y,\tau)}_-$. Thus $f^{(2)}$ can also be obtained and denoted as $V_4$. Consequently, the Lyapunov stability of each equilibrium can be derived from $V(\psi)=V_{\text{q}'}(\psi)+V_4(\psi)-\sum_{j=1}^n h(\psi_j(0))$.
 $\hfill\blacksquare$
\end{appendices}


\begin{thebibliography}{00}
\bibitem{freeman2001biocomplexity} W. J. Freeman and R. Kozma and P. J. Werbos,	\textit{Biocomplexity: adaptive behavior in complex stochastic dynamical systems}, Biosystems, 59(2001), pp. 109-123.
\bibitem{Tang2009} W. Ma, A. Trusina, H. El-Samad, et. al. 
\textit{Defining network topologies that can achieve biochemical adaptation}, Cell, 138(2009), pp. 760-773.
\bibitem{Feinberg2019}M. Feinberg, \textit{Foundations of chemical reaction network theory}, Springer, 2019.
\bibitem{hangos2018} L. Márton, G. Szederkényi, K. M. Hangos, \textit{Distributed control of interconnected Chemical Reaction Networks with delay,} Journal of Process Contr., 71(2018), pp: 52-62.

\bibitem{Ri2021} F.A. Rihan, \textit{Delay differential equations and applications to biology}, Singapore: Springer, 2021.
\bibitem{Li2023}Y. Yang, K. Foster, K. Z. Coyte, A. Li, \textit{Time delays modulate the stability of complex ecosystems}, Nat. Ecol. Evol., 7(2023), pp. 1610-1619.
\bibitem{Giordano2014} F. Blanchini, G. Giordano, \textit{Piecewise-linear Lyapunov functions for structural stability of biochemical networks}, Automatica, 50(2014), pp. 2482-2493.
\bibitem{Fang2019} Z. Fang, C. Gao. \textit{Lyapunov function partial differential equations for chemical reaction networks: Some special cases}, SIAM J. Appl. Dyn. Syst., 18(2019), pp. 1163-1199.
\bibitem{Fang2020} Z. Fang, A. Van Der Schaft, C. Gao, \textit{A graphic formulation of nonisothermal chemical reaction systems and the analysis of detailed balanced networks}, SIAM J. Appl. Dyn. Syst., 19(2020), pp. 2594-2627.
\bibitem{Horn1972} F. Horn, R. Jackson, \textit{General mass action
kinetics}, Arch. Ration. Mechan., 47(1972), pp: 81–116.
\bibitem{Feinberg1995} M. Feinberg, \textit{The existence and uniqueness of steady states for a class of chemical reaction networks}, Arch. Ration. Mechan., 132(1995), pp. 311-370.
\bibitem{Feinberg1972} M. Feinberg, \textit{Complex balancing in general kinetic systems}, Arch. Ration. Mechan., 49(1972), pp. 187-194.
\bibitem{Hangos2018} G. Lipták, K. M. Hangos, M. Pituk, \textit{Semistability of complex balanced kinetic systems with arbitrary time delays}, Syst. Control Lett., 114(2018), pp. 38-43.
\bibitem{Alonso2012} G. Szederkényi, J. R. Banga, A. A. Alonso, \textit{CRNreals: a toolbox for distinguishability and identifiability analysis of biochemical reaction networks}, Bioinformatics, 28(2012), pp. 1549-1550.
\bibitem{Ye2023}
Q. Ye, R. Xu, D. Li, et al, \textit{Integrating multi-modal deep learning on knowledge graph for the discovery of synergistic drug combinations against infectious diseases}, Cell Rep. Phys. Sci., 4(2023).
\bibitem{Siegel2011}
M. D. Johnston, D. Siegel, \textit{Linear conjugacy of chemical reaction networks}, J. Math. Chem., 49(2011), pp. 1263-1282.
\bibitem{Zhang2022} 
X. Zhang, C. Gao, D. Dochain, \textit{On stability of two kinds of delayed chemical reaction networks}, IFAC-PapersOnLine, 55(2022), pp. 14-20.
\bibitem{H2021}
H. Hong, J. Kim, M. A. Al-Radhawi, E. D. Sontag, et al, 
\textit{Derivation of stationary distributions
of biochemical reaction networks via structure transformation}, Commun. Biol., 4(2021), pp. 620.
\bibitem{Feinberg1977}
M. Feinberg, and F. J. Horn, \textit{Chemical mechanism structure and the coincidence of the stoichiometric and kinetic subspaces}, Arch. Ration. Mechan., 66(1977), pp. 83-97.
\bibitem{G2019}
G. Lipták, M. Pituk, K. M. Hangos, \textit{Modelling and stability analysis of complex balanced kinetic systems with distributed time delays}, J. Process Control, 84(2019), pp. 13-23.
\end{thebibliography}
\end{document}